\documentclass[11pt]{article}
\usepackage{graphicx}
\usepackage{latexsym}
\usepackage{amssymb}
\usepackage{amsfonts}
\usepackage{amsmath}
\usepackage{amsthm}

\newtheorem{theorem}{Theorem}[section]

\newtheorem{lemma}[theorem]{Lemma}
\newtheorem{prop}[theorem]{Proposition}

\theoremstyle{remark}
\newtheorem{remark}[theorem]{Remark}

\theoremstyle{remark}
\newtheorem{definition}[theorem]{Definition}

\theoremstyle{remark}
\newtheorem{example}[theorem]{Example}

\numberwithin{equation}{section}

\newcommand{\ep}{\epsilon}

\newcommand{\reals}{\mathbb{R}}

\newcommand{\nats}{\mathbb{N}}
\newcommand{\comps}{\mathbb{C}}

\newcommand{\disk}{{\mathbb{D}^2}}
\newcommand{\half}{{\mathbb{R}_{+}^2}}
\newcommand{\sphere}{{\mathbb{S}^2}}

\newcommand{\nbhd}{\mathcal{N}}
\newcommand{\cover}{\mathcal{B}}

\newcommand{\inv}{^{-1}}

\newcommand{\ovl}{\bar}

\newcommand{\til}{\widetilde}

\newcommand{\wh}{\widehat}

\newcommand{\Hdim}{\mathcal{H}}
\newcommand{\subeq}{\subseteq}
\newcommand{\supeq}{\supseteq}

\def\diam{\operatorname{diam}}
\def\card{\operatorname{card}}
\def\dist{\operatorname{dist}}
\def\im{\operatorname{im}}

\begin{document}
\title{Quasisymmetric parameterizations of two-dimensional \\ metric planes }
\author{Kevin Wildrick \thanks{The author is partially supported by ?}\\
University of Michigan Department of Mathematics\\
Ann Arbor, MI 48019\\ 
 kwildric@umich.edu}
\date{}
\maketitle

\begin{abstract}
The classical uniformization theorem states that any simply connected Riemann surface is conformally equivalent to the disk, the plane, or the sphere, each equipped with a standard conformal structure. We give a similar uniformization for Ahlfors 2-regular, linearly locally connected metric planes; instead of conformal equivalence, we are concerned with quasisymmetric equivalence.  
\end{abstract}

\setcounter{section}{0} 

\section{Introduction}
\indent

Quasisymmetric maps are a generalization of conformal maps of Euclidean space to the general metric space setting.  Analogous to the uniformization theorem for Riemann surfaces, the task of characterizing a given metric space up to quasisymmetry is of general interest.  The spaces $\reals^n$ and $\mathbb{S}^n$ equipped with the standard metric are of particular interest in this problem, partially because a self-homeomorphism of Euclidean space is quasiconformal if and only if it is quasisymmetric.  As a result, the theory of quasiconformal mappings provides a guiding light.  Tukia and V\"ais\"al\"a \cite{QSembed} gave a simple intrinsic characterization of metric spaces quasisymmetrically equivalent to $\mathbb{S}^1$: a metric space $X$ homeomorphic to $\mathbb{S}^1$ is quasisymmetrically equivalent to it if and only if $X$ is doubling and linearly locally connected ($LLC$).  A similar characterization exists for $\reals$.  For $n \geq 3$, a complete characterization of $\mathbb{S}^n$ and $\mathbb{R}^n$ has yet to be given, and examples of Semmes \cite{SemmesEx} have shown that the problem is exceedingly difficult.  

In this paper, we focus on the case $n=2$.  Bonk and Kleiner \cite{QSParam} found necessary and sufficient conditions for a metric space to be quasisymmetrically equivalent to $\sphere$.  Under the additional assumption of Ahlfors $2$-regularity, this characterization is the same as in the one-dimensional case.

 \begin{theorem}[Bonk, Kleiner] \label{B-K}  Let $X$ be an Ahlfors 2-regular metric space homeomorphic to $\sphere$.  Then $X$ is quasisymmetrically equivalent to $\sphere$ if and only if $X$ is linearly locally connected. \end{theorem}
 
The purpose of this paper is to extend this result to metric spaces homeomorphic to the plane.  Throughout, we will use $\sphere$, ${\sphere}^*$, $\disk$, $\reals^2$, and $\half$ to denote the sphere, the once-punctured sphere, the open unit disk, the plane, and the open half-plane respectively, each endowed with the metric inherited from the ambient Euclidean metric.  We will denote the completion of a metric space $X$ by $\ovl{X}$, and the metric boundary by $\partial{X}:=\ovl{X}-X$. Our main result is the following theorem.

\begin{theorem}\label{main} Let $X$ be an Ahlfors 2-regular and linearly locally connected metric space homeomorphic to the plane or the sphere.  
\begin{itemize}
\item[(i)] If $X$ is bounded and complete, then $X$ is quasisymmetrically equivalent to $\sphere$.
\item[(ii)] If $X$ is bounded and $\card(\partial{X})=1$, then $X$ is quasisymmetrically equivalent to ${\sphere}^*$.
\item[(iii)] If $X$ is bounded and $\card(\partial{X})\geq 2$, then $X$ is quasisymmetrically equivalent to $\disk$.
\item[(iv)] If $X$ is unbounded and complete, then $X$ is quasisymmetrically equivalent to $\reals^2$.
\item[(v)] If $X$ is unbounded and not complete, then $X$ is quasisymmetrically equivalent to $\half$. 
\end{itemize}

The statements (i),(ii),(iv), and (v) are quantitative in the sense that the distortion function of each quasisymmetry depends only on the constants associated to the Ahlfors 2-regularity and linear local connectedness conditions.   In statement (iii), the distortion function also depends on the ratio $\diam{X}/\diam{\partial{X}}$.  

Conversely, if $X$ is a metric space which is quasisymmetrically equivalent to any of $\sphere$, ${\sphere}^*$, $\disk$, $\reals^2$, and $\half$, then $X$ is linearly locally connected with constant depending only on the distortion function of the quasisymmetry.    
\end{theorem}

Theorem \ref{main} shows that in order to determine the quasisymmetry type of an Ahlfors 2-regular, linearly locally connected metric space homeomorphic to the plane, one need only know if $d$ is bounded, and (roughly) how many non-convergent Cauchy sequences exist.  As quasisymmetric homeomorphisms map bounded sets to bounded sets and Cauchy sequences to Cauchy sequences, this is in the minimal information required to make such a determination.  Example \ref{diam ratio} below  shows that the dependence of the distortion function of the quasisymmetry in Theorem \ref{main}(iii) on the ratio $\diam{\partial{X}}/\diam{X}$ cannot be avoided.  The final statement of Theorem \ref{main} is well-known and is discussed in remark \ref{converse} below.

Theorem \ref{B-K} has an interesting application to hyperbolic geometry.  A well-known conjecture of Cannon states that for every Gromov hyperbolic group $G$ with boundary at infinity $\partial_{\infty}G$ homeomorphic to $\sphere$, there exists a discrete, co-compact, and isometric action of $G$ on hyperbolic 3-space.  By a theorem of Sullivan \cite{Sullivan}, this conjecture is equivalent to the following statement: if $G$ is a Gromov hyperbolic group, then $\partial_{\infty}G$ is homeomorphic to $\sphere$ if and only if $\partial_{\infty}G$ is quasisymmetrically equivalent to $\sphere$. The boundary $\partial_{\infty}G$ of a Gromov hyperbolic group has a natural family of $LLC$ and Ahlfors regular metrics.  Thus, Theorem \ref{B-K} confirms Cannon's conjecture if one of these metrics is Ahlfors $2$-regular.  Since this is not always the case, it is of particular interest to relax the Ahlfors regularity assumptions in Theorems \ref{B-K} and \ref{main}.   Recent progress on this problem includes \cite{QMRigid} and \cite{ConfDim}.

Theorem \ref{B-K} is quantitative in the same sense as Theorem \ref{main}.  Theorem \ref{main}(i) is merely a rephrasing of Theorem \ref{B-K}, included for completeness of the statement.  The authors of \cite{QSParam} note that the methods used to prove Theorem \ref{B-K} can also be used to establish Theorem \ref{main}(iv).  However, this approach requires the use of technical tools such as $K$-approximations of metric spaces and a discrete modulus, and has not been carried out in detail.  The methods employed in this paper are substantially more elementary, provided that one accepts Theorem \ref{B-K}.  

An outline of the proof of Theorem \ref{main} is as follows.  Let $X$ be as in the hypotheses of Theorem \ref{main}, and suppose that $X$ is a bounded space.  Bounded and Ahlfors regular spaces are totally bounded.  Thus, if $X$ is complete, it is homeomorphic to $\sphere$. Theorem \ref{B-K} then applies, proving Theorem \ref{main} (i).  If  $\hbox{card}(\partial{X})=1$, then $X$ is homeomorphic to the plane. Furthermore, $\ovl{X}$ is homeomorphic to the one-point compactification of $X$, which is $\sphere$. Applying Theorem \ref{B-K} produces a quasisymmetric equivalence of $\ovl{X}$ and $\sphere$, which restricts to a quasisymmetric equivalence of $X$ and $\sphere^*$.  

If  $\hbox{card}(\partial{X})\geq 2$, we show that $\partial{X}$ is homeomorphic to a circle.  This step is the core of the paper, and is a consequence of the following more powerful theorem.

\begin{theorem}\label{Circle} Let $X$ be a $\lambda$-$LLC$ metric space homeomorphic to the disk.  If $\ovl{X}$ is compact and $\partial{X}$ contains at least two points, then $\partial{X}$ is homeomorphic to the circle $\mathbb{S}^1$ and is $\lambda'$-$LLC$, where $\lambda'$ depends only on $\lambda$.  If in addition, $\partial{X}$ is doubling, then it is quasisymmetrically equivalent to $\mathbb{S}^1$, and the distortion function of the quasisymmetry depends only on $\lambda$ and the doubling constant. \end{theorem}

Note that if $X$ is doubling, then $\partial{X}$ is doubling as well.  To prove Theorem \ref{Circle}, we study the delicate interaction between the topological and metric properties of $X$. We show that $\partial{X}$ is a locally connected metric continuum such that the removal of any one point does not disconnect the space, while the removal of any two points does disconnect the space.  A theorem of point-set topology states that such a space is homeomorphic to the circle \cite{Wilder}.  In fact, our proof is quantitative, which leads to the additional conclusions regarding the $LLC$ condition and quasisymmetry.
  
Once it is established that $\partial{X}$ is homeomorphic to the circle, we may isometrically embed $X$ into the ``doubled'' space  $X'$ which is obtained by gluing two copies of $\ovl{X}$ together along $\partial{X}$. The space $X'$ is homeomorphic to $\sphere$ and is again Ahlfors 2-regular and $LLC$, and so we may apply Theorem \ref{B-K} to it.  The image of $X$ under the resulting quasisymmetry is an $LLC$ domain in $\sphere$ with boundary homeomorphic to the circle.  It is well known that such a domain is quasisymmetrically equivalent to $\disk$ (see Theorem \ref{Gehring} below).  Composing the various quasisymmetries yields the desired result.

In the case that $X$ is unbounded, we construct a new metric on $X$ which results in a bounded metric space $\wh{X}$, which we call the ``warp" of $X$. This process, also employed in \cite{QMRigid}, is analogous to obtaining the standard (extrinsic) metric on $\sphere^*$ from the standard metric on $\reals^2$ via stereographic projection.  Similar warping processes for length spaces have recently been examined by Balogh and Buckley in \cite{warping}.  We show that $\wh{X}$ is again Ahlfors 2-regular and $LLC$, and that the boundary of $\wh{X}$ can be identified with $\partial{X} \cup \{\infty\}.$  Applying the bounded cases discussed above to $\wh{X}$ provides a quasisymmetry $\wh{f}\colon \wh{X} \to Y,$ where $Y$ is either $\sphere^*$ or $\disk$.  The warping process is designed so that the identity map $X \to \wh{X}$ is quasi-M\"obius.  This implies that  $\wh{f}$ descends to a quasisymmetry $f \colon X \to Z$, where $Z$ is $\reals^2$ or $\reals^2_+.$    

Theorem \ref{main}(iv) has already been used in \cite{Jacob} to relate the quasiconformal Jacobian problem to the classification of bi-Lipschitz images of the plane.  The quasiconformal Jacobian problem in the plane asks which non-negative locally integrable functions (weights) on $\reals^2$  are comparable to the Jacobian of a quasiconformal homeomorphism of $\reals^2$.  Suppose that given any weight, one could determine whether it is comparable to a quasiconformal Jacobian.  Then one can also determine whether a given metric space is bi-Lipschitz equivalent to the plane.  Let $(X,d)$ be a metric space.   If $X$ is bi-Lipschitz equivalent to the plane, then it is homeomorphic to the plane, Ahlfors $2$-regular, $LLC$, unbounded, and complete. Thus we may assume that $X$ satisfies the hypotheses of Theorem \ref{main}(iv).  Let $f \colon X \to \reals^2$ be the resulting quasisymmetric homeomorphism.  Elementary properties of quasisymmetric maps and a theorem of David and Semmes \cite{weights} show that the pushforward measure $\mu = f_*\Hdim^2$ satisfies 
$$d\mu(x)=w(x)dm_2(x)$$ 
for a so-called strong $A_{\infty}$-weight $w$ on $\reals^2$.  It is shown in \cite{Jacob} that $X$ is bi-Lipschitz equivalent to $\reals^2$ if and only if $w$ is comparable to the Jacobian of a quasiconformal homeomorphism of the plane. 

As mentioned above, Theorem \ref{main} can be viewed as a generalization of the classical uniformization theorem for Riemann surfaces.  One might also ask if other uniformization theorems can be similarly generalized.   It seems that techniques similar to those in this paper might be used to prove the following version of Koebe's uniformization on to circle domains.  Let $X$ be a bounded, Ahlfors $2$-regular, and $LLC$ metric space homeomorphic to a domain in $\sphere$ with $n$ boundary components.  Then $X$ is quasisymmetrically equivalent to $\sphere - \bigcup_{i=1}^n D_i$, where $\{D_i\}$ is a pairwise disjoint collection of closed balls or points.  In light of the work of He and Schramm \cite{Kobe}, one might also ask if such a theorem exists when countably many boundary components are allowed.  Bonk \cite{carpet} has recently given a result in this direction in the context of Sierpinski carpets.   

The techniques in this paper might also be used to prove a local version of Theorem \ref{main}.  Let $X$ be a proper and locally Ahlfors $2$-regular metric space which is homeomorphic to a surface.  Assume furthermore that $X$ is linearly locally contractible on compacta, i.e. that for every compact $K \subeq X$ there is a constant $\Lambda$ such that every ball $B(x, r)$ with $x \in K$ and $0 < r \leq \Lambda\inv$ is contractible inside of $B(x, \Lambda r).$  Then for each point $x \in X$, there is a neighborhood of $x$ which is quasisymmetrically equivalent to $\disk$. This statement plays a small role in the program of Heinonen et al.\cite{JuhaICM} to determine which submanifolds of $\reals^n$ are locally bi-Lipschitz equivalent to $\reals^2$.  

The author extends heartfelt thanks to his advisor Mario Bonk for his mentorship and many of the ideas in this paper.  Also, thanks to Juha Heinonen for many useful discussions.  

\section{Notation, Definitions, and Preliminary Results}
\indent

Where it will not cause confusion, we will refer to a metric space $(X,d)$ by $X$.  For $a \in X$ and $r>0$, we will use the following notation:
$$B_{X,d}(a,r):=\{x \in X\colon d(a,x)<r\},$$
$$\ovl{B}_{X,d}(a,r):=\{x \in X\colon d(a,x)\leq r\},$$
If $U\subeq X$ and $\ep>0$, we denote the $\epsilon$-neighborhood of $U$ in $X$ by $\nbhd_{\ep}^{X,d}(U)$.  We will often use $B(a,r)$, $B_d(a,r)$, or $B_X(a,r)$ in place of $B_{X,d}(a,r)$.  A similar convention will be used for closed balls, neighborhoods, and other objects which depend implicitly on the space $(X,d)$.  

Let $(X,d_X)$ and $(Y,d_Y)$ be metric spaces.  A homeomorphism $f\colon X \to Y$ is called quasisymmetric if there exists a homeomorphism $\eta\colon [0,\infty) \to [0,\infty)$ such that for all triples $a,b,c \in X$ of distinct points, we have 
$$\frac{d_Y(f(a),f(b))}{d_{Y}(f(a),f(c))} \leq \eta\left(\frac{d_X(a,b)}{d_X(a,c)}\right).$$
We will call the function $\eta$ the distortion function of $f$; when $\eta$ needs to be emphasized, we say that $f$ is $\eta$-quasisymmetric.  If $f$ is a quasisymmetric homeomorphism, then $f\inv$ is as well.  Thus we say that metric spaces $X$ and $Y$ are quasisymmetric or quasisymmetrically equivalent if there is a quasisymmetric homeomorphism from $X$ to $Y$.   We summarize some basic properties of quasisymmetric mappings in the following proposition.  Proofs can be found in \cite{QSembed} and \cite[Ch. 10]{LAMS}. 

\begin{prop}\label{QS properties} Let $f\colon X \to Y$ be an $\eta$-quasisymmetric homeomorphism of metric spaces. 
\begin{itemize}
\item[(i)] If $g\colon Y \to Z$ is a $\theta$-quasisymmetric homeomorphism, then $g \circ f$ is a $\theta \circ \eta$-quasisymmetric homeomorphism. 
\item[(ii)]  If $A \subeq B \subeq X$ are subsets with $0 < \diam{A} \leq \diam{B} < \infty$, then $\diam{f(B)}$ is finite and 
$$\frac{1}{2\eta\left(\frac{\diam{B}}{\diam{A}}\right)} \leq \frac{\diam{f(A)}}{\diam{f(B)}} \leq \eta\left(\frac{2\diam{A}}{\diam{B}}\right).$$
\item[(iii)] The map $f$ sends Cauchy sequences to Cauchy sequences, and there is a unique extension of $f$ to an $\eta$-quasisymmetric homeomorphism $\ovl{f}\colon \ovl{X} \to \ovl{Y}.$
\end{itemize}
\end{prop}

A homeomorphism of metric spaces $f\colon X \to Y$ is called quasi-M\"obius if there is a homeomorphism $\theta\colon [0,\infty) \to [0,\infty)$ such that for all quadruples $x_1,x_2,x_3,x_4 \in X$, the following relationship holds
$$[f(x_1),f(x_2),f(x_3),f(x_4)] \leq \theta([x_1,x_2,x_3,x_4]),$$
where the cross ratio is denoted
$$[x,y,z,w] := \frac{d(x,z)d(y,w)}{d(x,w)d(y,z)}.$$  
We will use the same notational conventions for quasi-M\"obius maps as for quasisymmetric maps.  The inverse of a quasi-M\"obius homeomorphism is again quasi-M\"obius, and a M\"obius transformation of $\reals^n$ is $\theta$-quasi-M\"obius with $\theta(t)=t$.  We will need the following result due to V\"ais\"al\"a \cite[Theorems 3.2 and 3.10]{QM}.

\begin{theorem}[V\"ais\"al\"a]\label{QM-QS}Let $f\colon X \to Y$ be a homeomorphism of metric spaces.  If $f$ is $\eta$-quasisymmetric, then $f$ is $\theta$-quasi-M\"obius with $\theta$ depending only on $\eta$.  If $X$ is unbounded and $f$ is a $\theta$-quasi-M\"obius homeomorphism which maps unbounded sequences to unbounded sequences, then $f$ is $\theta$-quasisymmetric. 
\end{theorem}

Let $\lambda > 1$.  A metric space $(X,d)$ is $\lambda$-linearly locally connected ($\lambda$-$LLC$) if for all $a \in X$ and  $r >0$ the following conditions are satisfied:
\begin{itemize}
\item[$(\lambda$-$LLC_1)$]  For each pair of distinct points $x,y \in B(a,r)$, there is a continuum $E \subeq B(a,\lambda r)$ such that $x,y\in E$,  
\item[$(\lambda$-$LLC_2)$] For each pair of distinct points $x,y \in X-B(a,r)$, there is a continuum $E \subeq X-B(a, r/\lambda)$ such that $x,y \in E$. 
\end{itemize}

Recall that a continuum is a connected, compact set containing more than one point.  Note that we do not place any upper restriction on the radius $r$ in this definition, though the $\lambda$-$LLC_2$ condition is vacuously true for $r> \diam(X,d)$.

\begin{remark}\label{LLC term}  The terminology ``linearly locally connected'' is justified by the following observation.  Suppose that $(X,d)$ is a $\lambda$-$LLC$ metric space, $x \in X$, and $r>0$.  Let $C(x)$ be the connected component of $B(x,r)$ containing $x$.  Then $B(x,r/\lambda)\subeq C(x) \subeq B(x,r)$.  
\end{remark}

V\"ais\"al\"a proved in \cite[Theorems 4.4 and 4.5]{QM} that the $LLC$ condition is preserved by quasi-M\"obius homeomorphisms.  In light of Theorem \ref{QM-QS}, we may state the following theorem. 

\begin{theorem}[V\"ais\"al\"a]\label{LLC is qs invariant} If $X$ is a $\lambda$-$LLC$ metric space and If $f\colon X \to Y$ is an $\eta$-quasisymmetric or $\eta$-quasi-M\"obius homeomorphism,  then $Y$ is $\lambda'$-$LLC$ for some $\lambda'$ depending only on $\lambda$ and $\eta$. 
\end{theorem}

\begin{remark}\label{converse} Each of the spaces $\sphere$, ${\sphere}^*$, $\disk$, $\reals^2$, and $\half$ is $LLC$.  This along with Theorem \ref{LLC is qs invariant} proves the final statement of theorem \ref{main}. 
\end{remark}

The question of which planar domains are quasisymmetrically equivalent to $\disk$ was essentially answered by Beurling and Ahlfors in \cite{Ahlfors-Beurling}.  However, the result was stated in terms of the $LLC$ condition by Gehring in \cite{Gehring}. 

\begin{theorem}[Beurling, Ahlfors, Gehring] \label{Gehring} Let $D \subeq \sphere$ be a domain which is $LLC$ when endowed with the standard metric, and such that $\partial{D}$ is connected and contains at least two points. Then there exists a quasisymmetric homeomorphism $f\colon D \to \disk$ with distortion function depending only on the $LLC$ constant of $D$.  
\end{theorem}

Let $I$ be any connected subset of $\reals$.  For any subset $U \subeq \ovl{X}$, we call a continuous map $\gamma\colon I \to U$ a path in $U$.  If the path $\gamma$ happens to be an embedding, then we call the image of $\gamma$ an arc in $U$.  We will make repeated use of the fact that the image of any path is arc-connected. A path $\gamma$ is called proper if for any compact set $K \subeq U$, the pre-image $\gamma\inv(K)$ is compact. The image of a path $\gamma$ will be denoted by $\im{\gamma}.$  

If $X$ is locally path-connected, we will often employ a condition similar to $LLC$ which uses arcs instead of continua.  This condition extends to the completion $\ovl{X}$ in a particularly nice way.  We say that a locally compact metric space $(X,d)$ is $\lambda$-$\til{LLC}$ if for all $a \in \ovl{X}$ and $r > 0$ the following conditions are satisfied:  
\begin{itemize} 
\item[$(\lambda$-$\til{LLC_1})$] For each pair of distinct points $x,y \in B_{\ovl{X}}(a,r)$, there is an embedding $\gamma\colon [0,1] \to \ovl{X}$ such that
$\gamma(0)=x$, $\gamma(1)=y$, $\gamma|_{(0,1)} \subeq X$, and $\gamma \subeq B_{\ovl{X}}(a,\lambda r),$  
\item[$(\lambda$-$\til{LLC_2})$] For each pair of distinct points  $x,y \in \ovl{X}-B_{\ovl{X}}(a,r)$, there is an embedding $\gamma\colon [0,1] \to \ovl{X}$ such that
$\gamma(0)=x$, $\gamma(1)=y$, $\gamma|_{(0,1)} \subeq X$, and $\gamma \subeq \ovl{X} - B_{\ovl{X}}(a,r/\lambda).$
\end{itemize}

If a metric space $X$ is $\lambda$-$\til{LLC}$, then it is also $\lambda$-$LLC$.  The next proposition states that the two conditions are quantitatively equivalent for the spaces in consideration in this paper.  

\begin{prop}\label{Better LLC}  Let $(X,d)$ be a locally compact, locally path-connected, and $\lambda$-$LLC$ metric space.   Then $X$ is $\lambda'$-$\til{LLC}$, where $\lambda'$ depends only on $\lambda$.  In particular, the space $\ovl{X}$ is $\lambda'$-$LLC$.  
\end{prop}
 
\begin{proof}  The key ingredient is the following statement:  If $U\subeq X$ is an open subset of $X$, and $E \subeq U$ is a continuum, then any pair of points $x,y \in E$ are contained in an arc in  $U$.  The details are straightforward and left to the reader. 
\end{proof}

The $\til{LLC}$ condition allows a useful addition to Remark \ref{LLC term}.

\begin{lemma}\label{closed-connected}    Let $(X,d)$ be a $\lambda$-$\til{LLC}$ metric space,  $p \in \partial{X}$, and $\ep > 0$.  Then there is a connected subset $C \subeq X$ which is closed in $X$, such that 
$$B_{\ovl{X}}(p, \ep/\lambda)\cap X \subeq C \subeq \ovl{B}_{\ovl{X}}(p,\ep)\cap X.$$
\end{lemma}
\begin{proof} 
Define 
$$S = \{(x,y) \in X \times X : x\neq y \ \hbox{and}\ x,y \in B_{\ovl{X}}(p,\ep/\lambda)\} \quad \hbox{and}\quad C_0= \bigcup_{(x,y) \in S} \gamma_{x,y},$$
where $\gamma_{x,y}$ is the arc connecting $x$ to $y$ provided by the $\lambda$-$\til{LLC}$ condition.
Taking $C$ to be the closure of $C_0$ in $X$ proves the lemma.
\end{proof}

Let $(X,d)$ be a metric space.  For any $Q \geq 0$, we define the $Q$-Hausdorff measure of a subset $E \subeq X$ by 
$$\Hdim^Q_d(E) := \lim_{\epsilon \to 0} \Hdim^{Q,\epsilon}_d(E)$$
where $\Hdim^{Q,\epsilon}_d(E)$ is the Carath\'{e}odory pre-measure defined as follows.  Let $\mathcal{B}_{\epsilon}$ be the collection of all covers $\mathcal{C}$ of $E$ by closed balls of radius less than $\epsilon$.  Then 
$$\Hdim^{Q,\epsilon}_d(E) := \inf_{\mathcal{C} \in \mathcal{B}_{\epsilon}} \sum_{B \in \mathcal{C}}(\text{radius}(B))^Q.$$
When computing $\Hdim^{Q,\epsilon}_d(E)$ it suffices to consider covers of $E$ by balls centered in the $\epsilon$-neighbor\-hood of $E$.  For a full description of Hausdorff measure and the Carath\'eodory construction, see \cite[Ch. 2.10]{Federer}.  Note that our definition differs from that in literature as we sum radii of balls rather than diameters of arbitrary closed sets; the resulting measures are comparable and thus equivalent for our purposes.  

A metric space $(X,d)$ is called Ahlfors $Q$-regular, $Q\geq 0$, if there exists a constant $K\geq 1$ such that for all $a \in X$ and $0<r\leq \diam{X}$, we have
\begin{equation}\label{strong reg def} \frac{r^Q}{K} \leq \Hdim^Q(\ovl{B}_d(a,r)) \leq Kr^Q.\end{equation}
\begin{remark}\label{ours is better} 
This is the definition used by Semmes in \cite{Semmes}, except that we do not require $X$ to be complete.  In \cite{QSParam}, Bonk and Kleiner use the slightly weaker condition that for all $a \in X$ and $0<r\leq \diam{X},$
\begin{equation}\label{weak reg def}  \frac{r^Q}{K} \leq \Hdim^Q({B}_d(a,r)) \leq Kr^Q .\end{equation}
The condition \eqref{strong reg def} implies the condition \eqref{weak reg def} with the same constant. In the case that $X$ is unbounded, the two conditions are equivalent.  
The main reason to use \eqref{strong reg def} rather than \eqref{weak reg def} is that \eqref{strong reg def} implies that 
$$\Hdim^Q(X) \leq K(\diam{X})^Q,$$
and so even for $r > \diam{X}$ we have the upper bound
$$\Hdim^Q(\ovl{B}_d(a,r)) \leq Kr^Q.$$
This is not true for spaces which only satisfy the weaker condition \eqref{weak reg def}. \end{remark}
  
A metric space $(X,d)$ is called $M$-doubling if for every $a \in X$ and all $r>0$, the open ball $B(a,r)$ can be covered by at most $M$ balls of radius $r/2$.  Note that there is no upper restriction on $r$ in this definition.   The next proposition lists some useful properties of Ahlfors $Q$-regular spaces. 

\begin{prop}\label{Qreg Prop} Let $Q \geq 0$, and let $(X,d)$ be an Ahlfors $Q$-regular metric space with constant $K$.  Then the following statements hold.
\begin{enumerate} 
\item[(i)] X is $M$-doubling where $M$ depends only on $K$ and $Q$.  
\item[(ii)] Any bounded subset of $X$ is totally bounded. 
\item[(iii)] The completion $\ovl{X}$ is Ahlfors $Q$-regular with constant $K'$ depending only on $K$ and $Q$.
\end{enumerate}
\end{prop}

In the proof of Proposition \ref{Qreg Prop}, we will need the following covering lemma, which is proven and discussed in \cite[Ch. 1]{LAMS}. 

\begin{lemma}[Basic Covering Lemma]\label{Basic Covering Lemma}
Let $(X,d)$ be a metric space.  Suppose that $\{B(x_i,r_i)\}_{i \in I}$ is a collection of balls in $X$ of uniformly bounded radius.  Then there exists a subset $J \subeq I$ such that 
$$\bigcup_{i \in I}B\left(x_i,r_i \right) \subeq \bigcup_{i \in J}B(x_i,5r_i),$$
and 
$$B\left(x_i,r_i\right) \cap B\left(x_j,r_j\right)= \emptyset$$
for distinct indices $i$ and $j$ in $J$. 
\end{lemma}

\begin{proof}[Proof of Proposition \ref{Qreg Prop}]  The proof of statement (i) is well-known and can be found, in particular, in \cite[Ch. 2.2]{Semmes}.   Statement (ii) follows directly from statement (i).  We will prove (iii).  Let $a \in \ovl{X}$ and let $r \leq \diam{\ovl{X}}=\diam{X}.$ We first consider the case that $a \in X$.  Let $\ep>0$ and consider any cover $\cover = \{\ovl{B}_{\ovl{X}}(x_i,r_i)\}_{i \in I}$ of $\ovl{B}_{\ovl{X}}(a,r)$ such that $r_i < \ep$.  If the center $x_i$ of a covering ball happens to be in the boundary $\partial{X}$, let $x_i'$ be any point in $\ovl{B}_{\ovl{X}}(x_i,r_i) \cap X$.  For those $x_i$ which are not boundary points, let $x_i'=x_i$.  Then for all $i \in I$, we have 
$$\ovl{B}_{\ovl{X}}(x_i,r_i) \subeq \ovl{B}_{\ovl{X}}(x_i',2r_i).$$
As a result, the collection 
$$\cover'=\{\ovl{B}_{\ovl{X}}(x_i',2r_i)\}_{i \in I}$$
is a cover of $\ovl{B}_{\ovl{X}}(a,r)$ by balls centered in $X$ of radius less than $2\ep$.  To prove the lower bound, we note that this implies the collection $\{\ovl{B}_{X}(x_i', 2r_i)\}_{i \in I}$ is a cover of $\ovl{B}_{X}(a, r)$ by balls in $X$ of radius less than $2\ep$.   As a result, 
$$\Hdim_{X}^{Q,2\ep}(\ovl{B}_{X}(a,r)) \leq \sum_{i \in I} (2r_i)^Q \leq 2^Q \sum_{i \in I} r_i^Q.$$
Since the cover $\mathcal{B}$ was arbitrary for the purposes of calculating $\Hdim_{\ovl{X}}^{Q, \ep}(\ovl{B}_{\ovl{X}}(a, r))$, letting $\ep$ tend to zero yields
$$\Hdim_X^Q(\ovl{B}_{X}(a, r)) \leq 2^Q \Hdim_{\ovl{X}}^Q(\ovl{B}_{\ovl{X}}(a, r)).$$
Thus, the $Q$-regularity of $X$ implies 
\begin{equation}\label{Hmeas interior lower}
\frac{r^Q}{2^QK} \leq \Hdim^Q_{\ovl{X}}(\ovl{B}_{\ovl{X}}(a,r)).
\end{equation}

To show the upper bound, we apply the Basic Covering Lemma to the collection $\cover'$.  Let
$\{\ovl{B}_{\ovl{X}}(x_i',10r_i)\}_{i \in J}$
be the resulting cover of $\ovl{B}_{\ovl{X}}(a,r).$ Now  
\begin{equation}\label{Hmeas interior upper 1}
\Hdim_{\ovl{X}}^{Q,10\ep}(\ovl{B}_{\ovl{X}}(a,r))\leq \sum_{i \in J} (10r_i)^Q \leq 5^Q \sum_{i \in J} \left(2r_i\right)^Q.
\end{equation}
For sufficiently small values of $\epsilon$, we have $\ovl{B}_X\left(x_i',2r_i\right)\subeq \ovl{B}_X(a,2r)$ for each $i \in J$.  Thus by \eqref{Hmeas interior upper 1}, the Ahlfors $Q$-regularity of $X$, and the disjointedness provided by the covering lemma, we have
$$\Hdim_{\ovl{X}}^{Q,10\ep}(\ovl{B}_{\ovl{X}}(a,r))\leq 5^Q \sum_{i \in J} K\Hdim^{Q}_{X}\left(\ovl{B}_{X}\left(x_i', 2r_i \right)\right) \leq 5^QK \Hdim^{Q}\left(\ovl{B}_{X}(a,2r)\right) \leq 10^Q K^2 r^Q.$$ 
Letting $\epsilon$ tend to zero yields that 
\begin{equation}\label{Hmeas interior upper}  \Hdim_{\ovl{X}}^{Q}(\ovl{B}_{\ovl{X}}(a,r)) \leq 10^Q K^2 r^Q.\end{equation}
Note that \eqref{Hmeas interior upper} holds if $r > \diam{X}$ as well.  

If $a \in \partial{X}$, we may find $a' \in X$ such that $d(a, a') < r/2.$  Then
$$\ovl{B}_{\ovl{X}}\left(a',\frac{r}{2}\right) \subeq \ovl{B}_{\ovl{X}}(a, r) \subeq \ovl{B}_{\ovl{X}}(a', 2r).$$ Then, using \eqref{Hmeas interior lower} and \eqref{Hmeas interior upper}, we have
\begin{equation}\label{Hmeas boundary}
 \frac{r^Q}{4^Q K} \leq  \Hdim^Q_{\ovl{X}}(\ovl{B}_{\ovl{X}}(a',r/2)) \leq \Hdim^Q_{\ovl{X}}(\ovl{B}_{\ovl{X}}(a,r)) \leq \Hdim^Q_{\ovl{X}}(\ovl{B}_{\ovl{X}}(a', 2r) \leq 20^Q K^2 r^Q.
\end{equation}
The estimates \eqref{Hmeas interior lower}, \eqref{Hmeas interior upper}, and \eqref{Hmeas boundary} show that $\ovl{X}$ is Ahlfors $Q$-regular with constant $20^QK^2$. \end{proof}

Given a $\lambda$-$LLC$ metric space $(X,d)$ which is Ahlfors $Q$-regular with constant $K$, we define the data of $(X,d)$ to be the triple $(\lambda, K, Q)$.  Propositions \ref{Qreg Prop} and \ref{Better LLC} show that if $(X,d)$ is an Ahlfors $Q$-regular and $LLC$ metric space, then $\ovl{X}$ is also Ahlfors $Q$-regular and $LLC$ with data depending only on the data of $X$.  

\section{The Sphere and Punctured Sphere}

\begin{proof}[Proof of Theorem \ref{main}(i)]  This is merely a restatement of Theorem \ref{B-K}.  Suppose that $X$ is a complete, bounded, Ahlfors $2$-regular, and $LLC$ metric space homeomorphic to the sphere or the plane.  By Proposition \ref{Qreg Prop} (ii), $X$ is totally bounded, and thus compact.  Accordingly, $X$ is homeomorphic to $\sphere$ and satisfies the hypotheses of Theorem \ref{B-K}, which provides an $\eta$-quasisymmetric homeomorphism  $f\colon X \to \sphere$ where $\eta$ depends only on the data of $X$.   
\end{proof}

\begin{proof}[Proof of Theorem \ref{main} (ii)] Suppose that $X$ is a bounded, Ahlfors $2$-regular, $LLC$ metric space such that $\card{\partial{X}}=1$.  By Proposition \ref{Qreg Prop} (ii) and (iii), and Proposition \ref{Better LLC},  $\ovl{X}$ is a compact, Ahlfors $2$-regular, and $LLC$ metric space with data depending only on the data of $X$.  A standard theorem of point-set topology \cite[Theorem 29.1]{Munkres} implies that $\ovl{X}$ is homeomorphic to the one-point compactification of $X$, namely $\sphere$.  Theorem \ref{B-K} now implies that there is a $\eta$-quasisymmetric homeomorphism $f\colon \ovl{X} \to \sphere$ where $\eta$ depends only on the data of $X$.  The restriction $f|_{X}$ is an $\eta$-quasisymmetric homeomorphism from $X$ to $\sphere^*.$      
\end{proof}

\section{The Boundary of a Disk}

By Proposition \ref{QS properties}(ii), a necessary condition for a metric space $X$ to be quasisymmetrically equivalent to $\disk$ is that $\partial{X}$ is a quasicircle, i.e. the quasisymmetric image of the circle $\mathbb{S}^1$. In this section, we prove Theorem \ref{Circle}, which provides sufficient conditions on $X$ for the boundary $\partial{X}$ to be a quasicircle.  

\begin{remark}\label{compact} Let $X$ be as in the assumptions of Theorem \ref{main}(iii); i.e. $X$ is homeomorphic to the plane, Ahlfors $2$-regular, $LLC$, bounded, and satisfies $\card{\partial{X}} \geq 2.$   By Proposition \ref{Qreg Prop}, the completion $\ovl{X}$ is doubling.  As a result, $\ovl{X}$ is compact and $\partial{X}$ is doubling.  Thus Theorem \ref{Circle} allows us to conclude that $\partial{X}$ is a quasicircle.  For the proof of Theorem \ref{main}(iii), we will only need the weaker conclusion that $\partial{X}$ is homeomorphic to the circle. 
\end{remark}

As mentioned in the introduction, we will show that $\partial{X}$ is homeomorphic to the circle by showing it is a locally connected metric continuum such that the removal of any one point does not disconnect the space, while the removal of any two points does disconnect the space.  We begin by giving the purely topological results which will be used in the proof.  

\begin{definition} Let $X$ be a topological space, and $k \in \nats$.  A subset $U \subeq X$ is $k$-ended if for every compact subset $K \subeq X$, there is another compact subset $K'\supeq K$ such that $U - K'$ has exactly $k$ components.  If, in addition, each of these components are arc-connected, then we say that $U$ is arc-$k$-ended. 
\end{definition}

\begin{remark}\label{compacta don't separate}  A trivial but useful example is the following: if $X$ is a topological space homeomorphic to the disk and $C\subeq X$ is a compact subset, then $X-C$ is arc-$1$-ended. 
\end{remark}

\begin{lemma}\label{proper} Let $X$ be a topological space homeomorphic to the disk, and suppose that $\gamma\colon (0,1) \to X$ is a proper embedding.  Then $X-\im\gamma$ has exactly two components, each of which is arc-connected and arc-$1$-ended.  Furthermore, there exists an ascending sequence $K_1 \subeq K_2 \subeq \hdots$ of compact subsets of $X$ with $\bigcup_{n \in \nats}K_n = X$ such that for each component $U$ of $X-\im\gamma$ and each $n \in \nats$, $U-K_n$ is arc-connected.
\end{lemma}

\begin{proof} As this is a purely topological result, we may assume that $X$ is $\sphere^*$, with the puncture at a point labeled $\infty$.  The assumption that $\gamma$ is proper now means that $\im\gamma \cup \{\infty\}$ defines a Jordan curve in $\sphere$. By the Sch\"onflies Theorem, there is a homeomorphism $\Theta\colon \sphere \to \sphere$ mapping $\im\gamma\cup\{\infty\}$ to a great circle $C$. Thus, up to homeomorphism, $X-\im\gamma$ is the complement of a line in $\reals^2$. In this case, the assertions of the lemma are clear. \end{proof}

\begin{lemma}\label{diff components} Let $X$ be a topological space homeomorphic to the disk, and let $\gamma$ and $\gamma'$ be proper embeddings of $(0,1)$ into $X$.  Suppose that there is a compact interval $I \subeq (0,1)$ such that $\gamma(t) = \gamma'(t)$ for all $t \in (0,1)-I.$  Then there is a compact subset $K$ of $X$ such that if $p, q \in X-K$ are in different components of $X-\im{\gamma'}$, then they are in different components of $X - \im{\gamma}.$  
\end{lemma}

\begin{proof} Let $U$ be a component of $X- \im{\gamma}.$  By Lemma \ref{proper}, we may find a compact set $K$ such that $\gamma(I) \subeq K$ and $U-K$ is arc-connected.  It suffices to show that if $p$ and $q$ are points in $U-K$, then they may be connected in $X- \im{\gamma'}.$  By assumption, $p$ and $q$ may be joined by an arc $\alpha$ which meets neither $K$ nor $\im{\gamma}$.  This implies that $\alpha$ does not meet $\im{\gamma'}$ either, and so $p$ and $q$ are in the same component of $X - \im{\gamma'}.$  
\end{proof}

A proof of the following statement may be found in \cite[V.2]{Newman}. 
\begin{lemma}\label{compact-closed}  Let $X$ be a topological space homeomorphic to the disk. Suppose that $K \subeq X$ is compact and connected, $C \subeq X$ is closed and connected, and $C \cap K$ is connected.  If $x$ and $y$ are points in $X$ which lie in the same component of $X - C$ and in the same component of $X - K$, then $x$ and $y$ also lie in the same component of $X - (C\cup K)$.  
\end{lemma}
  
We now turn to the proof of Theorem \ref{Circle}.  The $LLC$ condition is not needed to show that the boundary is a continuum. 

\begin{prop}\label{continuum} If $X$ is a metric space homeomorphic to the disk such that $\ovl{X}$ is compact, then $\partial{X}$ is a continuum if it contains at least two points.  \end{prop}

\begin{proof} As a closed subset of the compact space $\ovl{X}$, the boundary $\partial{X}$ is compact.  Assuming that $\partial{X}$ contains at least two points, it suffices to show that $\partial{X}$ is connected.  If $\partial{X}$ is not connected, we may find disjoint, non-empty, and closed subsets $A$ and $B$ of $\partial{X}$ with $A \cup B = \partial{X}$.  There is some $\ep>0$ such that $\dist(A,B)>2\ep.$  Let $U$ and $V$ be the $\ep$-neighborhoods of $A$ and $B$ respectively; then $U \cap V$ is empty.  Each point of $A$ is in the interior of $U$ by definition, and each point of $B$ is at a distance at least $\ep$ from $U$.  Thus $\partial{U} \cap \partial{X} = \emptyset$, and so $
\partial{U}$ is a compact subset of $X$. By Remark \ref{compacta don't separate}, there is a compact set $K \subeq X$ containing $\partial U$ such that each pair of points $u,v \in X - K$ can be connected by an arc which does not intersect $\partial U$.  Let $\delta = \dist(K,\partial{X})$, and let $\epsilon'<\min(\delta/2,\epsilon/2)$. We may find points $u$ and $v$ of $X$ in $\nbhd_{\ep'}(A)$ and $\nbhd_{\ep'}(B)$ respectively.  Then $u \in U$ and $v \in V$ but neither is in $K$. Thus they may be connected by an arc which does not intersect $\partial U$, contradicting the fact that $U$ and $V$ are disjoint.
\end{proof}

Throughout the rest of this section, we will assume that $X$ is a $\lambda$-$\til{LLC}$ metric space homeomorphic to the disk such that $\ovl{X}$ is compact and $\partial{X}$ contains at least two points.  Proposition \ref{Better LLC} shows that we have lost no generality in doing so.

\begin{prop}\label{minus 1} For each $p \in \partial{X}$, $\partial{X}-\{p\}$ is connected. 
\end{prop}

\begin{proof}  We argue by way of contradiction.  Suppose that $A$ and $B$ are disjoint, non-empty, and relatively closed subsets of $\partial{X}-\{p\}$ satisfying $A \cup B = \partial{X}-\{p\}$.  We may find disjoint open sets $U,V \subeq \ovl{X}$ containing $A$ and $B$ respectively.  Choose $\epsilon_1>0$ so small that we may find points $a \in A-B_{\ovl{X}}(p,\epsilon_1)$ and $b \in B-B_{\ovl{X}}(p,\epsilon_1).$  Let $$\ep_2=\frac{\ep_1}{2\lambda(4\lambda+1)}.$$  Then $X - (U \cup V \cup B_{\ovl{X}}(p,\epsilon_2/\lambda))$ is a compact subset of $X$.  Because $X$ is homeomorphic to the plane, there is a topological closed disk $K_1 \subeq X$ such that 
$$K_1 \supeq X - (U \cup V \cup B_{\ovl{X}}(p,\epsilon_2/\lambda)).$$
Note that $K_1$ is a compact and  connected subset of $X$.

By Lemma \ref{closed-connected}, there is a closed and connected subset $C_1 \subeq X$ such that 
$$B_{\ovl{X}}(p, \ep_2/\lambda)\cap X \subeq C_1 \subeq \ovl{B}_{\ovl{X}}(p,\ep_2)\cap X.$$  
We would like to apply Lemma \ref{compact-closed} to $C_1$ and $K_1$, but it may be
the case that $C_1 \cap K_1$ is not connected.  If $C_1 \cap K_1$ is connected, set $C=C_1$ and $K=K_1$.  Otherwise, we will add ``connectors'' to each component of $C_1\cap K_1$.  Set $\delta=\dist(C_1\cap K_1,\partial{X})$; note that $\delta \leq \ep_2$.  By compactness there is a cover of $C_1 \cap K_1$ by a finite collection of balls $\{B_i\}$ where $B_i:=B_X(x_i,\delta/2\lambda)$ with $x_i \in C_1 \cap K_1$.  Let
$C(x_i)$ denote the closure in $X$ of the component of $x_i$ in $B_X(x_i,\delta/2)$.  By Remark \ref{LLC term}, $B_i \subeq C(x_i)$, and so the collection
$\{C(x_i)\}$ is a cover of $C_1 \cap K_1$ by finitely many connected sets in $X$.  For each pair of distinct indices $i,j$, the $\lambda$-$\til{LLC}$ property provides an arc $\gamma_{ij}$ in $X$ connecting $x_i$ to $x_j$ inside $B_X(x_i,2\lambda d(x_i,x_j))$.  Since $C_1 \subeq \ovl{B}_{\ovl{X}}(p,\ep_2),$ we have 
$$\gamma_{ij} \subeq B_X(x_i,2\lambda d(x_i,x_j)) \subeq B_X(x_i,4\lambda\ep_2) \subeq B_{\ovl{X}}(p,\ep_1/2\lambda)\cap X.$$
We also have 
$$ \bigcup_i C(x_i) \subeq \bigcup_{i} B_X(x_i, \ep_2/2) \subeq B_{\ovl{X}}(p, \ep_1/2\lambda) \cap X.$$
Let 
$$K=K_1 \cup \bigcup_i C(x_i) \cup \bigcup_{i\neq j} \gamma_{ij} \quad\hbox{and}\quad C=C_1 \cup \bigcup_i C(x_i) \cup \bigcup_{i\neq j} \gamma_{ij}.$$ 
Now we have that $K$ is compact and connected, $C$ is closed and connected, and $C \cap K$ is connected. As $K$ is compact and the points $a$ and $b$ are in the boundary $\partial{X}$, Remark \ref{compacta don't separate} implies that we may find points $u$ and $v$ in the same component of $X-K$ such that 
$$u \in U \cap (X - B_{\ovl{X}}(p,\ep_1/2)) \quad \hbox{and} \quad v \in V \cap (X - B_{\ovl{X}}(p,\ep_1/2))$$
Furthermore, $C \subeq B_{\ovl{X}}(p,\ep_1/2\lambda)\cap X$, and so by $\lambda$-$\til{LLC}_2$ we see that $u$ and $v$ are in the same component of $X-C$.  Therefore, Lemma \ref{compact-closed} implies that $u$ and $v$ are in the same component of $X-(C\cup K)$.  However, 
$$C \cup K \supeq  \left(B_{\ovl{X}}(p,\ep_2/\lambda)\cap X\right) \cup \left(X-(U \cup V \cup B_{\ovl{X}}(p,\ep_2/\lambda)\right)\supeq X - (U \cup V).$$
This means that $u$ and $v$ lie in a connected subset of $U \cup V$, which contradicts the facts that $u \in U$, $v \in
V$, and $U \cap V = \emptyset$.  Thus $\partial{X}-\{p\}$ is connected. 
\end{proof}

\begin{definition}\label{crosscut}  Let $p$ and $q$ be distinct points in $\partial{X}$.  A crosscut connecting $p$ and $q$ is an embedding $\gamma\colon [0,1] \to X \cup \{p,q\}$ such that $\gamma(0)=p$ and $\gamma(1)=q$. \end{definition}

Note that if $\gamma$ is a crosscut, then $\gamma:(0,1) \to X$ is a proper embedding.  

\begin{lemma}\label{big proper arc lemma}   Let $\gamma$ be any crosscut, and let $U$ and $V$ be the components of $X - \im{\gamma}.$ The following statements hold:
\begin{itemize}
\item[(i)] $\ovl{X}=\ovl{U} \cup \ovl{V},$
\item[(ii)] $\ovl{U}-\im{\gamma}$ and $\ovl{V}-\im{\gamma}$ are the components of $\ovl{X}-\im{\gamma}$,  
\item[(iii)] $\ovl{U}\cap \partial{X}$ and $\ovl{V} \cap \partial{X}$ are connected.
\end{itemize}
\end{lemma}

\begin{proof} (i)  This follows immediately from the definitions.  

(ii)  To show that $\ovl{U}-\im{\gamma}$ and $\ovl{V}-\im{\gamma}$ are the components of $\ovl{X}-\im{\gamma}$, we assert that they are each relatively closed and connected, they do not intersect, and their union is all of $\ovl{X}-\im{\gamma}.$  Clearly they are relatively closed, and by (i) we have $\ovl{U} \cup \ovl{V} = \ovl{X}$.

By definition, $U$ is connected and does not intersect $\im{\gamma}.$   Thus, to show that $\ovl{U}-\im{\gamma}$ is connected, it suffices to show that each point $u \in \ovl{U}-\im{\gamma}$ may be connected to some point $u' \in U$ without crossing $\im{\gamma}$.  Let $u \in \ovl{X}-\im{\gamma}.$ Since $\im{\gamma}$ is a compact subset of $\ovl{X}$, there is an $\ep>0$ such that $B_{\ovl{X}}(u, \ep)$ does not meet $\im{\gamma}$.  There is a point $u' \in U$ with $d(u,u') < \ep/\lambda$.  The $\lambda$-$\til{LLC}$ condition provides an arc $\alpha$ connecting $u$ to $u'$ inside of $B_{\ovl{X}}(u, \ep)$, and so $\alpha$ does not intersect $\im{\gamma}$.   The same argument shows that $\ovl{V} - \im{\gamma}$ is connected, and a slightly modified version shows that $\ovl{U} \cap \ovl{V} - \im{\gamma}$ is empty.

(ii) We now show that $\ovl{U} \cap \partial{X}$ is connected; the same argument will apply to $\ovl{V} \cap \partial{X}$. Suppose that $C$ and $D$ are disjoint, non-empty, and closed subsets of $\ovl{U} \cap \partial{X}$ such that $C \cup D = \ovl{U} \cap \partial{X}$.  Then we may find an $\epsilon>0$ such that $\dist(C,D)>\ep$. 

We claim that there is some $\delta>0$ such that 
$$U \cap \nbhd_{\delta}(\partial{X}) \subeq \nbhd_{\ep}(C)\cup\nbhd_{\ep}(D).$$
If not, then for all $n \in \nats$, we may find points $u_n \in U$ and $x_n \in \partial{X}$ such that 
$$d(u_n,x_n)<\frac{1}{n} \quad \hbox{and} \quad \dist(x_n, C \cup D)\geq \dist(u_n, C \cup D) - d(u_n, x_n) \geq \ep-\frac{1}{n}.$$
As $\partial{X}$ is compact, there is a subsequence $\{x_{n_k}\}$ converging to a point $x \in \partial X - (\nbhd_{\ep}(C)\cup\nbhd_{\ep}(D)).$  By construction $\{u_{n_k}\}$ also converges to $x$, contradicting the fact that $\ovl{U} \cap \partial{X} = C \cup D,$ and the claim is proven.  

Consider that $K:=X-\nbhd_{\delta}(\partial{X})$ is a compact subset of $X$.  By Lemma \ref{proper}, $U$ is one-ended, and so we may find a compact subset $K' \supeq K$ such that $U - K'$ is connected.  However, the claim shows that $\nbhd_{\ep}(C)$ and $\nbhd_{\ep}(D)$ constitute a cover of $U-K'$ by non-empty, disjoint, open sets.  This is a contradiction, and so $\ovl{U}\cap \partial{X}$ must be connected. 
\end{proof}

\begin{prop}\label{minus 2}  For any pair of distinct points $p,q \in \partial{X}$, the set $\partial{X}-\{p,q\}$ is not connected.  
\end{prop}

\begin{proof} The $\lambda$-$\til{LLC}$ condition provides a crosscut $\gamma$ connecting $p$ to $q$. Let $U$ and $V$ be the components of $X-\im{\gamma},$  and let $A=\ovl{U} \cap \partial{X}-\{p,q\}$ and $B = \ovl{V} \cap \partial{X}-\{p,q\}$.  Then $A$ and $B$ are relatively closed subsets of $\partial{X}-\{p, q\}$.  It follows from \ref{big proper arc lemma}(i) that $A \cup B = \partial{X}- \{p,q\}$.  We will show $A \cap B = \emptyset$ and that $A$ and $B$ are non-empty.  

Suppose that there is some point $z \in A \cap B$. Then $z \notin \im\gamma$, and so there is some $\epsilon >0$ such that $B_{\ovl{X}}(z,\epsilon)\cap \im\gamma =\emptyset.$ By assumption, we may find points $u \in U$ and $v \in V$ contained in $B_{\ovl{X}}(z, \epsilon/\lambda)$.  By the $\lambda$-$\til{LLC}$ condition, there is an arc connecting $u$ to $v$ which does not intersect $\im{\gamma}$.  This is a contradiction, and so $A \cap B = \emptyset$.

By symmetry, it suffices to show that $A$ is non-empty.  Lemma \ref{proper} provides an exhaustion $\{K_n\}_{n \in \nats}$ of $X$ by compact sets such that for each $n \in \nats$, $U-K_n$ are is arc-connected.  For each $n \in \nats$, we may find points 
$$p_n \in \im\gamma \cap (X-K_n) \cap B_{\ovl{X}}(p,1/n) \quad \hbox{and} \quad q_n \in \im\gamma \cap (X-K_n) \cap B_{\ovl{X}}(q,1/n).$$   
Because each $K_n$ is compact, there is a sequence of positive numbers $\{\epsilon_n\}_{n \in \nats}$ tending to zero such that $B_X(p_n,\epsilon_n)
\cap K_n = \emptyset$ and $B_X(q_n, \ep_n)\cap K_n=\emptyset$.  Each point on $\im{\gamma}$ is a limit point of  $U$, so there are points $u_n \in B_X(p_n,\epsilon_n) \cap U$ and $u_n' \in B_X(q_n,\ep_n)\cap U$. We may connect $u_n$ to $u_n'$ via an arc $\gamma_n$ in $U - K_n$.   By the connectedness of $\gamma_n$, there is a point $x_n \in \gamma_n \cap U$ such that 
$$d(u_n,x_n)=\frac{d(u_n,u'_n)}{2}\leq d(u'_n,x_n).$$
By the compactness of $\ovl{X}$, $\{x_n\}_{n \in \nats}$ subconverges to some point $x \in \ovl{X}$.   Since $\bigcup_{n \in \nats} K_n = X$, we have $x \in \partial{X}$.  Furthermore, $u_n \to p$ and $u_n' \to q$ as $n \to\infty$, so
 $$d(p,x)=\frac{d(p,q)}{2}\leq d(q,x).$$
Thus $x \in \partial{X}-\{p,q\}$.  Since $\{x_n\}\subeq U$, we have $x \in A$. 
 \end{proof}
 
\begin{lemma}\label{agree off compact}  Let $a, b, p, q$ be distinct points in $\partial{X}$, and let $\gamma$ and $\gamma'$ be crosscuts connecting $a$ and $b$.  Suppose that there is a compact interval $I \subeq (0,1)$ such that $\gamma(t)=\gamma'(t)$ for all $t \in [0,1]-I$.  If $p$ and $q$ are in different components of $\ovl{X}-\im{\gamma'}$, then they are in different components of $\ovl{X}-\im{\gamma}$. 
\end{lemma}

\begin{proof}
 By Lemma \ref{diff components}, there is a compact set $K$ such that if $p'$ and  $q'$ are points of $X-K$ which are in different components of $X-\im{\gamma'}$, then they are in different components of $X- \im{\gamma}.$  We may find $\ep>0$ such that 
$$B_{\ovl{X}}(p, \ep) \cap \left(\im{\gamma'} \cup \im{\gamma} \cup K \right) = \emptyset $$  
Let $p' \in X$ such that $d(p,p') < \ep/\lambda.$  Then $p' \notin K$.  The $\lambda$-$\til{LLC}$ condition provides and arc $\alpha$ connecting $p$ to $p'$ which does not intersect $\im{\gamma} \cup \im{\gamma'}.$  
Thus $p$ and $p'$ are in the same component of $\ovl{X} - \im{\gamma'}$ and the same component of $\ovl{X} - \im{\gamma}$.  Similarly, we may find $q' \in X-K$ such that $q$ and $q'$ are in the same component of $\ovl{X} - \im{\gamma'}$ and the same component of $\ovl{X} - \im{\gamma}.$  Thus, if $p$ and $q$ are in different components of $\ovl{X}-\gamma'$, then $p'$ and $q'$ are in different components of $\ovl{X}-\im{\gamma'}$.  This implies that $p'$ and $q'$ are in different components of $X-\im{\gamma'}$, and so by Lemma \ref{diff components} they are in different components of $X-\im{\gamma}.$  Lemma \ref{big proper arc lemma}(ii)  allows us to conclude that $p'$ and $q'$ are in different components of $\ovl{X}-\im{\gamma}$, and hence so are $p$ and $q$.    
\end{proof}

\begin{lemma}\label{4 points}  Let $a,b,p,q$ be distinct points in $\partial{X}$. Let 
$\gamma_{pq}$ be a crosscut connecting $p$ to $q$ and $\gamma_{ab}$ be a crosscut connecting $a$ to $b$.  If $a$ and $b$ are in different components of $\ovl{X}-\im{\gamma_{pq}}$, then $p$ and $q$ are in different components of $\ovl{X}-\im{\gamma_{ab}}$. \end{lemma}

\begin{proof}  The fact that the points $a,b,p,q$ are all distinct implies that $K = \gamma_{a,b}\cap \gamma_{p,q}$ is a compact subset of $X$. For ease of notation, we identify $\reals^2$ with $\comps$ in this proof.  Let $\phi\colon X \to \comps$ be any homeomorphism.  We may find $R>0$ such that $\phi(K) \subeq B_{\comps}(0, R).$   As $\gamma_{ab}|_{(0,1)}$ and $\gamma_{pq}|_{(0,1)}$ are proper, we may find parameters $t_1,t_2,s_1,s_2 \in (0,1)$ such that 
$$t_1 = \min\{t\in (0,1): |\phi \circ \gamma_{ab}(t)|=R\} \quad \hbox{and}\quad t_2 = \max\{t\in (0,1):|\phi\circ \gamma_{ab}(t)|=R\},$$
$$s_1 = \min\{s\in (0,1):|\phi \circ \gamma_{pq}(t)|=R\} \quad \hbox{and}\quad s_2 = \max\{s\in (0,1):|\phi \circ \gamma_{pq}(t)|=R\}.$$
Set 
$$a' = \phi \circ \gamma_{ab}(t_1) = Re^{i\theta_{a'}} \quad \hbox{and} \quad b'=\phi\circ \gamma_{ab}(t_2) = Re^{i\theta_{b'}},$$ 
$$p' = \phi\circ\gamma_{pq}(s_1)= Re^{i\theta_{p'}} \quad \hbox{and} \quad q'=\phi\circ\gamma_{pq}(s_2)= Re^{i\theta_{q'}},$$
where $\theta_{a'},\theta_{b'},\theta_{p'},\theta_{q'} \in [0, 2\pi).$  Note that $a',b', p', q'$ must be distinct points on the circle $\mathcal{C}=\{z \in \reals^2: |z| = R\}$, for otherwise $\im{\gamma_{ab}}$ and $\im{\gamma_{pq}}$ are disjoint, contradicting the assumption that $a$ and $b$ are in different components of $\ovl{X}-\im{\gamma_{pq}}$.   Thus we may consider the cyclic order of $a', b', p', q'$ on the circle $C$. Suppose that $\theta_{a'} < \theta_{b'} < \theta_{p'} < \theta_{q'}.$ Define arcs $\alpha$ and $\beta$ in $X$ by 
$$\alpha = \phi\inv(\{Re^{i\theta} \in \comps: \theta_{a'} \leq \theta \leq \theta_{b'} \} )\quad \hbox{and} \quad \beta =  \phi\inv(\{Re^{i\theta} \in \comps: \theta_{p'} \leq \theta \leq \theta_{q'} \}),$$
and set 
$$\gamma'_{ab} = \im(\gamma_{ab}|_{[0,t_1)}) \cup \alpha \cup \im(\gamma_{ab}|_{(t_2,1]}) \quad \hbox{and} \quad \gamma'_{pq} = \im(\gamma_{pq}|_{[0,s_1)}) \cup \beta \cup \im(\gamma_{pq}|_{(s_2,1]}).$$
Then $\gamma'_{ab}$ connects $a$ to $b$ without intersecting $\gamma'_{pq}$.  However, Lemma \ref{agree off compact} implies that the points $a$ and $b$ are in different components of $\ovl{X}-\gamma'_{pq},$ yielding a contradiction.  Any case where $a'$ and $b'$ are adjacent in the cyclic order on $C$ can be handled in a similar fashion.

Now suppose that $\theta_{a'} < \theta_{p'} < \theta_{b'} < \theta_{q'}.$  In this case, define $\alpha \subeq X$ to be the inverse image under $\phi$ of the line segment in $\comps$ from $a'$ to $b'$.  Similarly, let $\beta \subeq X$ be the inverse image under $\phi$ of the line segment from $p'$ to $q'$.   Define $\gamma'_{ab}$  and $\gamma'_{pq}$ as before.  Then it is clear that $p'$ and $q'$ are in different components of $\ovl{X}-\im{\gamma'_{ab}}$.  Furthermore, $p'$ can be connected to $p$ without intersecting $\im{\gamma'_{ab}}$, and similarly for $q'$ and $q$.  Thus $p$ and $q$ are in different components of $\ovl{X}-\im{\gamma'_{ab}}.$  Lemma \ref{agree off compact} then implies that they are in different components of $\ovl{X}-\im{\gamma_{ab}}$.  Any case where $a'$ and $b'$ are not adjacent in the cyclic order on $C$ can be handled in a similar fashion.  \end{proof}

\begin{prop}\label{locally connected}
The boundary $\partial{X}$ is $LLC_1$ with constant depending only on the $LLC$ constant of $X$.
\end{prop}

\begin{proof} Let $p \in \partial{X}$, and $r>0$.  It suffices to find a continuum $E$ such that 
$$B_{\ovl{X}}(p, r) \cap \partial{X} \subeq E  \subeq B_{\ovl{X}}(p, 4\lambda^4 r) \cap \partial{X}.$$

By Proposition \ref{continuum}, we may assume that there is some point $q\in \partial{X}-B_{\ovl{X}}(p, 4\lambda^4 r)$.  The $\til{LLC}$ condition provides a crosscut $\gamma_{pq}$ connecting $p$ to $q$.  Let $U$ and $V$ be the components of $X-\im{\gamma_{pq}}$, and set $A=\ovl{U}\cap \partial{X}$ and $B=\ovl{V}\cap \partial{X}$.  By Lemma \ref{big proper arc lemma}(iii), $A$ and $B$ are connected.  As $\{p,q\} = A\cap B$ and $d(p,q) > 4\lambda^4 r$, we may find distinct points $a  \in A$ and $b \in B$ such that $d(p,a)= 2\lambda^2 r$ and  $d(p,b) =  2\lambda^2 r$.    By Lemma \ref{big proper arc lemma}(ii), the points $a$ and $b$ lie in different components of $\ovl{X}-\im{\gamma_{pq}}.$  The $\lambda$-$\til{LLC}_1$ condition provides a crosscut $\gamma_{ab}$ connecting $a$ to $b$ with $\im{\gamma_{ab}} \subeq B_{\ovl{X}}(p,3\lambda^3 r).$  Let $W$ be the component of $X-\im{\gamma_{ab}}$ with $p \in \ovl{W}\cap \partial{X}$.  Set $E:=\ovl{W}\cap \partial{X}$.  Applying Lemma \ref{big proper arc lemma} again,  we see that $E$ is a continuum.  We first show that $E \subeq B_{\ovl{X}}(p, 4\lambda^4 r) \cap \partial{X}.$  Suppose that there is a point $x \in E - B_{\ovl{X}}(p, 4\lambda^4 r)$.  By the $\lambda$-$\til{LLC}_2$ condition, there is a path connecting $x$ to $q$ without intersecting $B_{\ovl{X}}(p, 4\lambda^3 r)$.  This implies that $p$ and $q$ are in the same component of $\ovl{X}-\gamma_{ab}$, contradicting Lemma \ref{4 points}. 

We now show that $B_{\ovl{X}}(p, r) \cap \partial{X} \subeq E.$  Since $\gamma_{ab}$ is uniformly continuous, we may find parameters $0<t_a<1$ and $0< t_b < 1$ such that 
$$\diam(\gamma_{ab}([0,t_a]))\leq \lambda^2 r \quad \hbox{and} \quad \diam(\gamma_{ab}([t_b, 1]))\leq \lambda^2 r.$$
Set $a' = \gamma_{ab}(t_a)$ and $b' =\gamma_{ab}(t_b)$.  Then $a',b' \in X - B_{\ovl{X}}(p, \lambda^2 r),$ and so the $\lambda$-$\til{LLC}_2$ condition provides an embedding $\gamma_{a',b'}\colon [0,1] \to X$ such that $\gamma(0)=a'$, $\gamma(1)=b'$, and $\im{\gamma} \subeq X-B_{\ovl{X}}(p, \lambda r).$  Consider that the set 
$$S = \gamma_{ab}([0, t_a]) \cup \im{\gamma_{a'b'}} \cup \gamma_{ab}([t_b,1])$$
does not intersect $B_{\ovl{X}}(p, \lambda r),$ and is the image of a path in $\ovl{X}$.  Since the image of a path in $\ovl{X}$ is arc-connected, we may find a crosscut $\gamma'$ connecting $a$ to $b$ with $\im{\gamma'} \subeq S$.   Furthermore, after re-parameterization, we may find a compact interval $I \subeq (0,1)$ such that $\gamma'(t)=\gamma_{ab}(t)$ for all $t \in [0,1] - I$.   Suppose that there is a point $x \in B_{\ovl{X}}(p, r) \cap \partial{X}$ which is not contained in $E$.  Then $x$ and $p$ are in different components of $\ovl{X}-\im{\gamma_{ab}}.$ By Lemma \ref{agree off compact}, this implies that $x$ and $p$ are in different components $\ovl{X}-\im{\gamma'}$. However, the $\lambda$-$\til{LLC}_1$ condition shows that $x$ and $p$ may be connected by an arc contained in $B_{\ovl{X}}(p,\lambda r)$.  This is a contradiction.
\end{proof}

\begin{remark}\label{LLC2 from LLC1} If $Y$ is a $\lambda$-$LLC_1$ metric space homeomorphic to the circle $\mathbb{S}^1$,  then $Y$ is also $\lambda'$-$LLC$, where $\lambda'$ depends only on $\lambda$.  This follows quickly from the following elementary topological facts about the circle $\mathbb{S}^1$:
\begin{itemize}
\item If $x,y,z,w$ are distinct points in $\mathbb{S}^1$, then $x$ and $y$ are in different components of $\mathbb{S}^1- \{z, w\}$ if and only if $z$ and $w$ are in different components of $\mathbb{S}^1-\{x,y\}$.  
\item If $E \subeq \mathbb{S}^1$ is a continuum containing points $x, y \in \mathbb{S}^1$, then $E$ contains at least one of the components of $\mathbb{S}^1-\{x,y\}$.  
\end{itemize}
\end{remark}

\begin{proof}[Proof of Theorem \ref{Circle}] Let $X$ be a bounded, $\lambda$-$LLC$ metric space homeomorphic to the disk.  By Propositions \ref{continuum}, \ref{minus 1}, and \ref{minus 2}, the boundary $\partial{X}$ is a metric continuum such the removal of one point does not separate the space, while the removal of two points does separate the space.  Proposition \ref{locally connected} implies that, in particular, $\partial{X}$ is locally connected.  A recognition theorem of point-set topology \cite{Wilder} states that such a space is homeomorphic to the circle $\mathbb{S}^1$.   Proposition \ref{locally connected} and Remark \ref{LLC2 from LLC1} show that $\partial{X}$ is $\lambda'$-$LLC$ where $\lambda'$ depends only on $\lambda$. 

Tukia and V\"ais\"al\"a \cite{QSembed} characterized metric spaces quasisymmetrically equivalent to $\mathbb{S}^1$ in the following way: if $Y$ is a metric space homeomorphic to $\mathbb{S}^1$, then $Y$ is quasisymmetrically equivalent to $\mathbb{S}^1$ if and only if it is doubling and satisfies the $LLC_1$ condition.  Furthermore, they show that the distortion function of the quasisymmetry depends only on the $LLC_1$ and doubling constants.  This proves the final statement of Theorem \ref{Circle}.   
\end{proof}

\section{The Disk}

Throughout this section, let $(X,d)$ be a locally compact, bounded, and incomplete metric space.  Let $X'$ be the space obtained by gluing two copies of $\ovl{X}$ together by the identity map along $\partial{X}$.  We will denote elements of $X'$ by $[x,i]$ where $x \in X$ and $i \in \{1,2\}$; if $x \in \partial{X}$, then we will use the notation $[x,i]=[x]=[x,j].$  If $E \subeq \ovl{X}$, we set $[E,i]:=\{[x,i]: x \in E\}.$  By local compactness, we have $\dist(x, \partial{X})>0$ for each $x \in X$.  Thus there is a natural metric $d'$ on the space $X'$ given by 
$$d'([x,i],[y,j]):= \begin{cases}
			d(x,y) & i=j \\
			\inf\{d(x,z)+d(z,y): z \in \partial{X}\} & i\neq j \\
		    \end{cases}.$$
Note that $\diam{X'} \leq 2 \diam{\ovl{X}}$, and that $X$ embeds isometrically in $X'$.    
\begin{remark} \label{projection} The triangle inequality shows the projection map $[x,j] \mapsto x$ does not increase distance.  
\end{remark}

\begin{prop}\label{double}  Suppose that $(X,d)$ is an Ahlfors $Q$-regular and $LLC$.  Then the metric space $(X',d')$ is Ahlfors $Q$-regular and $LLC$, with data depending only on $Q$, the data of $(X,d)$, and the ratio $\diam{X}/\diam{\partial{X}}.$
\end{prop}

\begin{proof}   We begin by showing that $(X',d')$ is Ahlfors $Q$-regular.  By Proposition \ref{Qreg Prop}, we may assume that $\ovl{X}$ is Ahlfors $Q$-regular with constant $K$.  Let $[z,i] \in X'$, and let $r \leq \diam{X'}.$   We first give a lower estimate for $\Hdim^Q_{X'}(\ovl{B}_{X'}([z,i],r)).$  Note that $\frac{1}{2}r \leq \diam{\ovl{X}}$. Let $\ep > 0$, and consider any cover $\{\ovl{B}_{X'}([x_n,i_n],r_n)\} $ of $\ovl{B}_{X'}([z,i],r/2)$ by closed balls in $X'$ of radius less than $\ep$.  Then $\{\ovl{B}_{\ovl{X}}(x_n,r_n)\}$ is a cover of $\ovl{B}_{\ovl{X}}(z,r/2)$ by closed balls in $\ovl{X}$ of radius less than $\ep$.  Thus
 \begin{equation}\label{double reg lower} \frac{r^Q}{2^QK} \leq \Hdim^{Q}_{\ovl{X}}(\ovl{B}_{\ovl{X}}(z,r/2)) \leq  \Hdim^{Q}_{X'}(\ovl{B}_{X'}([z,i],r/2)) \leq \Hdim^{Q}_{X'}(\ovl{B}_{X'}([z,i],r)). \end{equation}

We now show an upper estimate.  If $\{\ovl{B}_{\ovl{X}}(x_n,r_n)\}$ is any cover of $\ovl{B}_{\ovl{X}}(z,r)$ by closed balls in $\ovl{X}$ of radius less than $\ep$, then $\{\ovl{B}_{X'}([x_n,1],r_n)\} \cup  \{\ovl{B}_{X'}([x_n,2],r_n)\}$ is a cover of $\ovl{B}_{X'}([z],r)$ by closed balls in $X'$ of radius less than $\ep$.  Therefore,
\begin{equation}\label{double reg upper 2}  \Hdim^{Q}_{X'}(\ovl{B}_{X'}([z],r)) \leq 2Kr^Q. \end{equation}
Combining \eqref{double reg lower} and \eqref{double reg upper 2},  we see that $X'$ is Ahlfors $Q$-regular with constant $\max(2^QK, 2K)$.

We now show that $(X',d')$ is $LLC$.  By Proposition \ref{Better LLC}, we may assume that $X$ is $\lambda$-$\til{LLC}$.   Let $[z,i] \in X'$ and $r >0 $.  Let $[x,j]$ and $[y,k]$ be points in  $B_{X'}([z,i],r)$.  By Remark \ref{projection} we have $x,y \in B_{\ovl{X}}(z,r)$.  

First suppose that $z \in X$ and $r < \dist(z, \partial{X}).$ This implies that $i=j=k$.  The $\lambda$-$\til{LLC}$ condition on $\ovl{X}$ provides a continuum $E_{xy} \subeq B_{\ovl{X}}(z,\lambda r)$ containing $x$ and $y$.     Then $[E_{xy},i] \subeq B_{X'}([z,i],\lambda r) $ is a continuum in $X'$ containing $[x,j]$ and $[y,k]$. 

Next, suppose that $z \in \partial{X}$ and $r >0$.  The $\lambda$-$\til{LLC}$ condition on $\ovl{X}$ provides continua $E_{xz}$ and $E_{yz}$ contained in $B_{\ovl{X}}(z, \lambda)$ and containing $\{x,z\}$ and $\{y,z\}$ respectively.  Now
$[E_{xz}, j] \cup [E_{yz},k] \subeq B_{X'}([z],\lambda r)$ is a continuum in $X'$ containing $[x,j]$ and $[y,k]$. 

Finally we consider the case that $z \in X$ and $r \geq \dist(z, \partial{X}).$  There is a point $z' \in \partial{X}$ such that $d(z,z') < 2r.$  Consider that 
$$B_{X'}([z,i],r) \subeq B_{X'}([z'], 3r) \subeq B_{X'}([z'], 3\lambda  r) \subeq B_{X'} ([z,i], (3\lambda  +2)r).$$
These inclusions and and the discussion above show that there is a continuum $E$  containing $[x,j]$ and $[y,k]$ such that $E \subeq B_{X'} ([z,i], (3\lambda  +2)r).$  We have now shown that $X'$ is $3\lambda +2$-$LLC_1$.    

Next, we show the $LLC_2$ condition.  Let $[x,j]$ and $[y,k]$ be points in $X' - B_{X'}([z,i],r).$  First suppose that $z \in X$ and $r < \dist(z,\partial{X})$.  This implies that neither $x$ nor $y$ are in $B_{\ovl{X}}(z,r).$  If $z'$ is any point of $\partial{X}$,  the $\lambda$-$\til{LLC}$ condition provides continua $E_{xz'}$ and $E_{yz'}$  contained in $\ovl{X}-B_{\ovl{X}}(z, r/\lambda )$ and containing $\{x,z'\}$ and $\{y,z'\}$ respectively .  Then $[E_{xz'},j] \cup [E_{yz'},k] \subeq X' - B_{X'}([z,i], r/\lambda )$ is a continuum in $X'$ containing $[x,j]$ and $[y,k]$.  
  
Now suppose that $z \in \partial{X}$.  We must have that $r \leq \diam{X'} \leq 2 \diam{X}$, for otherwise $X' - B_{X'}([z],r)=\emptyset$.  Setting $\alpha = \diam{X}/\diam{\partial{X}}$, we have 
$$\frac{r}{8\alpha} \leq \frac{\diam{\partial{X}}}{4} < \frac{\diam{\partial{X}}}{2}.$$ 
Thus there is a point $z' \in \partial{X} - B_{\ovl{X}}(z,  r/8\alpha).$  Furthermore, neither $x$ nor $y$ is in $B_{\ovl{X}}(z,r/8\alpha).$ Thus the $\lambda $-$LLC_2$ condition on $\ovl{X}$ provides continua $E_{xz'}$ and $E_{yz'}$  contained in $\ovl{X}-B_{\ovl{X}}(z, r/8\alpha \lambda )$ and containing ${x,z'}$ and ${y,z'}$ respectively .  Now, $[E_{xz'},j] \cup [E_{yz'},k] \subeq X' - B_{X'}([z,i], r/8\alpha \lambda )$ is a continuum in $X'$ containing $[x,j]$ and $[y,k]$.

We return to the case $z \in X$, and now allow that $r < 64\alpha \lambda \dist(z, \partial{X}).$  Then $[x,j]$ and $[y,k]$ are in the complement of $B_{X'}([z,i],r/64\alpha\lambda ).$  The first case above provides a continuum $E \subeq X' - B_{X'}([z,i],r/64\alpha \lambda ^2)$ containing $[x,j]$ and $[y,k]$.   If $r \geq 64\alpha \lambda \dist(z,\partial{X})$, then we may find a point $z' \in \partial{X}$ such that $d(z,z')< r/(32\alpha\lambda )$.  This implies that $[x,j]$ and $[y,k]$ are in the complement of $B_{X'}\left([z'],r/2\right).$  Now, the second case above provides a continuum $E$ containing $[x,j]$ and $[y,k]$ such that 
$$E \subeq X' - B_{X'}\left([z'], \frac{r}{16\alpha\lambda }\right) \subeq X'- B_{X'}\left([z,i],\frac{r}{32\alpha\lambda }\right).$$ 
We have now shown that $X'$ is $64\alpha\lambda^2$-$LLC_2$.    
\end{proof} 

\begin{example}\label{diam ratio}  The Ahlfors $2$-regularity and $LLC_1$ constants of $(X',d')$ do not depend on $\alpha$. However, the dependence of the $LLC_2$ constant  of $(X',d')$ on this ratio cannot be avoided.  Let $\epsilon >0$ and consider $X_{\ep}=\sphere-\ovl{B}(a,\epsilon).$  The $LLC_2$ constant of $X_{\ep}$ does not depend on $\epsilon$, while the $LLC_2$ constant of $X'_{\ep}$ tends to infinity as $\epsilon$ tends to zero.   The spaces $X_{\ep}$ also show that the distortion function of the uniformizing quasisymmetry provided by Theorem \ref{main} depends on $\alpha$ as well.  To see this, suppose that $f_{\ep}\colon X_{\ep} \to \disk$ is an $\eta$-quasisymmetric homeomorphism.  The map $f_{\ep}$ extends to an $\eta$-quasisymmetric homeomorphism $\ovl{f}_{\ep}\colon \ovl{X_{\ep}} \to \ovl{\mathbb{D}}^2$ sending $\partial{X_{\ep}}$ to $\partial{\disk}$.   For sufficiently small $\ep$, Proposition \ref{QS properties} shows that
$$1 = \frac{\diam{\partial{\disk}}}{\diam{\ovl{\mathbb{D}}^2}} \leq \eta\left(\frac{2\diam{\partial{X_{\ep}}}}{\diam{\ovl{X}_{\ep}}}\right) = \eta(2\ep).$$
Letting $\ep$ tend to zero yields a contradiction. 
\end{example}

\begin{proof}[Proof of  Theorem \ref{main} (iii).]  Suppose that $X$ is an Ahlfors $2$-regular, $LLC$, and bounded metric space homeomorphic to the plane, with $\card{\partial{X}}\geq 2$.  Remark \ref{compact} shows that we may apply Theorem \ref{Circle} and conclude that the boundary $\partial{X}$ is homeomorphic to the circle.  Theorem \ref{double} shows that $X'$ is an Ahlfors $2$-regular, $LLC$ metric space homeomorphic to $\sphere$ with data depending only on the data of $X$ and the ratio $\diam{X}/\diam{\partial{X}}$.  Theorem \ref{B-K} provides a quasisymmetric homeomorphism $f\colon X' \to \sphere$ whose distortion function depends only on the data of $X'$.  Theorem \ref{LLC is qs invariant} shows that $f(X)$ is an $LLC$ disk inside $\sphere$, and Proposition \ref{QS properties}(iii) shows that $\partial{f(X)}$ is connected and contains more than two points.  Theorem \ref{Gehring} provides a quasisymmetric homeomorphism $g\colon f(X) \to \disk$ whose distortion function depends only on the $LLC$-constant of $f(X)$, and hence only on the data of $X$ and the ratio $\diam{X}/\diam{\partial{X}}$.  The map $g\circ f$ is the desired quasisymmetric homeomorphism.    
\end{proof}

\section{The Plane and Half-Plane}

Throughout this section, let $(X,d)$ be a connected and unbounded metric space.  We wish to ``warp" $(X,d)$ to create a bounded metric space.   This warping process, which was also employed in \cite{QMRigid}, is analogous to obtaining the standard extrinsic metric on $\sphere$ from the standard metric on $\reals^2$.

Fix a basepoint $p \in X$, and define for all $x,y \in X$
$$\rho_p(x,y):= \frac{d(x,y)}{(1+d(x,p))(1+d(y,p))}.$$
In general, $\rho_p$ is not a metric on $X$.   To force the triangle inequality, we define 
$$\wh{d}_p(x,y) = \inf \sum_{i=0}^{k-1} \rho_p(x_i,x_{i+1}),$$
where the infimum is taken over all finite sequences of points $x=x_0,x_1,...,x_k=y$ in $X$.

As $p \in X$ will remain fixed throughout, we will suppress the reference to $p$ in the definitions above, using instead $\wh{d}=\wh{d}_p$ and $\rho=\rho_p$.  For further ease of notation, for all $x \in X$ we set 
$$h(x):=\frac{1}{1+d(x,p)}.$$

\begin{remark}\label{h ineq}
Note that for any $u,v \in X$ we have 
$$ \left| h(u)-h(v) \right| = h(u)h(v)\left| d(v,p)- d(u,p) \right| \leq \rho(u,v).$$
This and the triangle inequality show that for any sequence $x=x_0,x_1,\hdots,x_k=y$ of points in $X$, we have 
$$\sum_{i = 0}^{k-1}\rho(x_i,x_{i+1}) \geq |h(x)-h(y)|.$$
Thus for all points $x$ and $y$ in $X$ we have
$$\wh{d}(x,y) \geq |h(x)-h(y)|.$$
\end{remark} 

In order to show that $\wh{d}$ is a metric on $X$, we need the following lemma which is proven in \cite{QMRigid}. 
\begin{lemma}\label{1/4 ineq}  For all $x,y \in X$, we have
$$\frac{1}{4} \rho(x,y)\leq \wh{d}(x,y) \leq \rho(x,y).$$
\end{lemma}

Lemma \ref{1/4 ineq} shows that $\wh{d}(x,y)=0$ implies $x=y$ for all points $x$ and $y$ in $X$.  It follows from the definitions that $\wh{d}$ is symmetric and satisfies the triangle inequality.  Thus, we define the ``warped" version of $(X,d)$ to be the metric space $(X,\wh{d})$.   In this warped space, distances from $p$ may be calculated from the $d$-distance from $p$.

\begin{lemma}\label{p symmetry 1} If $x \in X$, then $\wh{d}(x,p)=1-h(x).$
\end{lemma}  

\begin{proof} From Lemma \ref{1/4 ineq}, we see that 
$\wh{d}(x,p) \leq \rho(x,p) = 1-h(x).$
On the other hand, setting $y=p$ in Remark \ref{h ineq} shows that 
$\wh{d}(x,p)\geq |h(x)-1|=1-h(x).$  \end{proof}

Let $\wh{X}$ denote the completion of $(X,\wh{d})$ and let $\wh{\partial}X=\wh{X}-X$ be the metric boundary of $(X,\wh{d})$. We seek a description of $\wh{\partial}X$ in terms of $\partial{X}$. 

\begin{lemma}\label{never complete}  Let $\{x_n\}$ be a sequence in $X$. If $d(x_n, p) \to \infty$, then $\{x_n\}$ is a non-convergent $\wh{d}$-Cauchy sequence.  Conversely, if $\{x_n\}$ is a $\wh{d}$-Cauchy sequence which is $d$-unbounded, then $d(x_n, p) \to \infty.$
\end{lemma}

\begin{proof} Suppose that $\{x_n\} \subeq X$ satisfies $d(x_n, p) \to \infty$.  Let $\epsilon >0$; we may find some integer $N > 0$ such that if $n \geq N$, then 
$$\frac{1}{1+d(x_n,p)} < \frac{\ep}{2}.$$ 
For $n \geq N$ and $k$ any positive integer,  Lemma \ref{1/4 ineq} and the triangle inequality show that 
$$\wh{d}(x_n,x_{n + k}) \leq \frac{d(x_n,p) + d(x_{n+k},p)}{(1+d(x_n,p))(1+d(x_{n+k},p))} \leq \frac{1}{1+d(x_{n+k},p)} + \frac{1}{1+d(x_n,p)} < \ep.$$
This shows that $\{x_n\}$ is a $\wh{d}$-Cauchy sequence.  Suppose there is some $x \in X$ such that $\wh{d}(x_n,x) \to 0$.  Then by Lemma \ref{p symmetry 1}
$$1-h(x) = \wh{d}(x,p)= \lim_{n \to \infty}\wh{d}(x_n,p)=1 - \lim_{n \to \infty}h(x_n) = 1.$$
This implies that $h(x)=0$, contradicting that $x \in X$. 
 
 Now, let $\{x_n\} \subeq X$ be a $\wh{d}$-Cauchy sequence which is $d$-unbounded.  If $d(x_n, p)$ does not tend to infinity, then there exists some $R \geq 0$ such that for infinitely many positive integers $n$ we have $d(x_n,p)<R$.   By Lemma \ref{p symmetry 1}
$$ \wh{d}(x_n,p) = \frac{d(x_n,p)}{1+d(x_n,p)}, $$
and so there are infinitely many positive integers $n$ such that $x_n \in B_{\wh{d}}(p,R/(1+R)).$
On the other hand, $\{x_n\}$ is $d$-unbounded, so there are infinitely many positive integers $n$ such that $d(x_n,p)>2R$.  As a result, there are infinitely many positive integers $n$ such that $x_n \notin B_{\wh{d}}(p,2R/(1+2R)).$  This contradicts the assumption that $\{x_n\}$ is a Cauchy sequence in $(X,\wh{d})$.   \end{proof}

If $\{x_n\}$ and $\{y_n\}$ are $\wh{d}$-Cauchy sequences which are $d$-unbounded, then by Lemma \ref{never complete} the sequence $\{x_1,y_1,x_2,y_2, \hdots \}$ has no $d$-bounded subsequences and is again a $\wh{d}$-Cauchy sequence.  Thus we may define a distinguished point $\infty \in \wh{\partial}X$ corresponding to $\wh{d}$-Cauchy sequences which are $d$-unbounded.  There is a special relationship between the basepoint $p$ and the point $\infty \in \wh{\partial}{X}.$

\begin{lemma} \label{p symmetry 2} If $x \in X$, then $\wh{d}(x,\infty)=h(x).$
\end{lemma}

\begin{proof} We must show that if $\{y_n\} \subeq X$ is any $\wh{d}$-Cauchy sequence which is $d$-unbounded, then $\wh{d}(x,y_n) \to h(x)$ as $n \to \infty$.  Lemma \ref{1/4 ineq} shows that 
$$\lim_{n \to \infty} \wh{d}(x,y_n) \leq \lim_{n \to \infty} \rho(x,y_n) = h(x).$$
On the other hand, by Remark \ref{h ineq} we have
$$\lim_{n \to \infty}\wh{d}(x,y_n) \geq \lim_{n \to \infty}|h(x)-h(y_n)|.$$
By Lemma \ref{never complete}, $\{y_n\}$ satisfies $d(y_n, p) \to \infty$, so $h(y_n)$ tends to zero.  Thus 
$$\lim_{n \to \infty}\wh{d}(x,y_n) \geq h(x).$$
\end{proof}

\begin{remark}\label{infty balls}Since $h(x)=\frac{1}{1 + d(x,p)}$, Lemmas \ref{p symmetry 1} and \ref{p symmetry 2} provide the first step in relating the metrics $d$ and $\wh{d}$.  An example of their usefulness is the equality
$$B_{\wh{d}}(\infty,r)=\begin{cases}
\wh{X} & r>1\\
\wh{X}-\{p\} & r = 1 \\
\wh{X}-\ovl{B}_d(p,\frac{1-r}{r}) & r < 1 \\
\end{cases}.$$
In particular, this shows that $\diam{(\wh{X},\wh{d})} \leq 2.$  It is also convenient to record that for all $r>0$,
$$X - B_{d}(p,r)=\ovl{B}_{\wh{d}}\left(\infty,\frac{1}{1+r}\right) \cap X.$$
\end{remark}

It is not possible to give such an exact description of every $\wh{d}$-ball, but the following lemma shows that the metrics $d$ and $\wh{d}$ are ``comparable away from infinity''.   This fact is the essential ingredient in showing that if $(X,d)$ is an Ahlfors $Q$-regular and $LLC$ space, then $(X,\wh{d})$ also has these properties.  

\begin{lemma} \label{inclusions}  Let $C>1$.  If $a \in X$, and $r \leq \frac{\wh{d}(a,\infty)}{C}$, then 
\begin{equation}\label{inclusions hat middle}B_{X,d}\left( a,\frac{r}{\wh{d}(a,\infty)^2}\cdot\frac{C}{C+1}\right) \subeq B_{X, \wh{d}}(a,r) \subeq B_{X,d}\left( a,\frac{r}{\wh{d}(a,\infty)^2}\cdot\frac{4C}{C-1}\right)\end{equation}
Furthermore, if $R \leq \frac{1}{\wh{d}(a,\infty)(C+1)}$, then
\begin{equation} \label{inclusions d middle} B_{X, \wh{d}}\left(a,R\wh{d}(a,\infty)^2\cdot\frac{C-1}{4C}\right) \subeq B_{X,d}(a,R) \subeq B_{X,\wh{d}}\left(a,R\wh{d}(a,\infty)^2\cdot\frac{C+1}{C}\right).\end{equation}
The inclusions \eqref{inclusions hat middle} and \eqref{inclusions d middle} also hold when all open balls are replaced with closed balls. 
\end{lemma}

\begin{proof} Suppose $x \in X$ satisfies $\wh{d}(a,x)\leq \frac{\wh{d}(a,\infty)}{C}$ for some $C>1$.  By Lemma \ref{p symmetry 2} we have
$$d(a,x) = \frac{\rho(a,x)}{\wh{d}(a,\infty)\wh{d}(x,\infty)}.$$
Lemma \ref{1/4 ineq} and the triangle inequality imply that 
$$\frac{\wh{d}(a,x)}{\wh{d}(a,\infty)( \wh{d}(a,\infty) + \wh{d}(a,x))} \leq d(a,x) \leq \frac{4\wh{d}(a,x)}{\wh{d}(a,\infty)( \wh{d}(a,\infty) - \wh{d}(a,x))}.$$
Since $\wh{d}(a,\infty)\geq C\wh{d}(a,x),$ this yields
$$\left(\frac{C}{C+1}\right)\frac{\wh{d}(a,x)}{\wh{d}(a,\infty)^2}\leq d(a,x) \leq \left(\frac{4C}{C-1}\right)\frac{\wh{d}(a,x)}{\wh{d}(a,\infty)^2},$$
and \eqref{inclusions hat middle} follows.  The inclusions \eqref{inclusions d middle} follow from \eqref{inclusions hat middle}  by noting that if $R \leq \frac{1}{\wh{d}(a,\infty)(C+1)},$ then 
$$\max\left(R\wh{d}(a,\infty)^2\cdot\frac{C-1}{4C}, R\wh{d}(a,\infty)^2\cdot\frac{C+1}{C}\right) \leq \frac{\wh{d}(a,\infty)}{C}.$$
\end{proof}

\begin{lemma}\label{pre hat boundary} Let $\{x_n\} \subeq X$ be a $d$-bounded sequence.  Then $\{x_n\}$ is a $d$-Cauchy sequence if and only if it is a $\wh{d}$-Cauchy sequence.  Furthermore $\{x_n\}$ $d$-converges to a point $x \in X$ if and only if it $\wh{d}$-converges to $x$. \end{lemma} 

\begin{proof}  By Lemma \ref{1/4 ineq}, we see that $\wh{d}(x,y) \leq d(x,y)$ for all points $x$ and $y$ in $X$.  This implies that if $\{x_n\}$ is a $d$-Cauchy sequence, then it is a $\wh{d}$-Cauchy sequence.  Furthermore,  this shows that if $\{x_n\}$ $d$-converges to a point $x \in X$, it $\wh{d}$-converges to $x$ as well.   

Now suppose that $\{x_n\}$ is a $\wh{d}$-Cauchy sequence which is $d$-bounded.   There is some $R>0$ such that $\{x_n\}\subeq B_{d}(p,R)$.   By Lemma \ref{p symmetry 2}, this implies that for all $n \in \nats$,  
\begin{equation}\label{hat boundary eq} \frac{1}{2(1+R)} < \frac{\wh{d}(x_n,\infty)}{2}. 
\end{equation}
Let $\ep >0$.  As $\{x_n\}$ is a $\wh{d}$-Cauchy sequence, there is some $N$ such that for all $n \geq N$ we have $x_n \in B_{\wh{d}}(x_N, \ep')$, where 
$$\ep' := \min\left( \frac{\ep}{8(1+R)^2} , \frac{1}{2(1+R)}\right).$$
Inequality \eqref{hat boundary eq} shows that we may apply the inclusion \eqref{inclusions hat middle} to $B_{\wh{d}}(x_N, \ep')$ with constant $C=2$.  Thus we see that for all $n \geq N$ 
$$x_n \in B_{d}\left(x_N, \frac{8\ep'}{\wh{d}(x_N,\infty)^2} \right) \subeq B_{d}(x_N,\ep).$$
This shows that $\{x_n\}$ is a $d$-Cauchy sequence and that if $\{x_n\}$ $\wh{d}$-converges to a point $x \in X$, it $d$-converges to $x$ as well. 
\end{proof}

\begin{prop}\label{hat boundary}  There is a bijection between $\wh{\partial}{X}$ and $\partial{X}\cup \{\infty\}$.   
\end{prop}

\begin{proof}  This follows from Lemmas \ref{never complete} and \ref{pre hat boundary}.
\end{proof}

\begin{lemma}\label{QM} The identity map $\iota\colon (X,d) \to (X,\wh{d})$ is a $\theta$-quasi-M\"obius homeomorphism with $\theta(t)=16t$.
\end{lemma} 

\begin{proof} That the identity map is a homeomorphism follows from Lemma \ref{pre hat boundary}; that it is $16t$-quasi-M\"obius follows from Lemma \ref{1/4 ineq}. 
\end{proof}

We now make rigorous the statement that the warping process is analogous to obtaining the extrinsic metric on $\sphere$ from the standard metric on $\reals^2$.  

\begin{lemma}\label{plane to sphere} Let $|\cdot |_{\reals^n}$ denote the standard metric structure on $\reals^n$, and $\wh{|\cdot |}_{\reals^n}$ denote the corresponding warped metric with basepoint at the origin.  Then $(\reals^2,\wh{|\cdot|}_{\reals^2})$ is bi-Lipschitz equivalent to the extrinsic metric on $\sphere^*.$    
\end{lemma}

\begin{proof}  Let $s\colon \reals^2 \to \sphere-\{0,0,1\}$ be the stereographic projection map. 
By Lemma \ref{1/4 ineq}, it suffices to show that there is a constant $L>1$ such that for all $x,y \in \reals^2$
$$\frac{1}{L}\frac{|x-y|_{\reals^2}}{(1+|x|)(1+|y|)} \leq |s(x)-s(y)|_{\reals^3} \leq L \frac{|x-y|_{\reals^2}}{(1+|x|)(1+|y|)}.$$
A calculation shows that 
$$|s(x)-s(y)|_{\reals^3} = \frac{2|x-y|_{\reals^2}}{\sqrt{(1+|x|^2)(1+|y|^2)}},$$ 
and the result follows with $L=4.$  \end{proof}

We now have the tools needed to prove that the warping procedure preserves Ahlfors $Q$-regularity and the $LLC$ condition quantitatively.  

\begin{prop}\label{hat regular}  Let $(X,d)$ be a connected and unbounded metric space, and let $Q>0$.  If $(X,d)$ is Ahlfors Q-regular with constant $K$, then $(X,\wh{d})$ is Ahlfors Q-regular with a constant depending only on $Q$ and $K$. 
\end{prop}

\begin{proof}  Throughout this proof, we will only consider balls centered in $X$ as objects in $X$, not in the completion $\ovl{X}$.  Thus we will use the notation 
$$\ovl{B}_{d}(a,r)=\ovl{B}_{X,d}(a,r) \quad \hbox{and} \quad \ovl{B}_{\wh{d}}(a,r) = \ovl{B}_{X,\wh{d}}(a,r).$$

The general method is to construct a cover of a ball in one metric from a cover of a ball in the other metric, in a quantitative way.  The main tool is Lemma \ref{inclusions}.   For technical reasons which will become clear later, we fix any $\til{C}>1$ such that 
$$C := \frac{8\til{C}(\til{C}+1)}{\til{C}-1} >2.$$  
Let $a \in X$ and $r \leq \diam{(X,\wh{d})}$.  The first step is to estimate $\Hdim^Q_{\wh{d}}(\ovl{B}_{\wh{d}}(a,r))$ in the case  that $r \leq \frac{\wh{d}{(a,\infty)}}{C}.$   Since $C>2$, we may fix $\ep>0$ such that 
$$\epsilon < \min\left(\frac{\wh{d}(a,\infty)-2r}{C}, r\right).$$
Suppose that $\ovl{B}_{\wh{d}}(a,r)$  is covered by a collection of closed balls $\{\wh{B}_i\}_{i \in I}$ where $\wh{B}_i := \ovl{B}_{\wh{d}}(x_i, r_i)$, $r_i < \epsilon,$ and $\wh{d}(a, x_i) < r + \epsilon$ for each $i \in I$.  From this cover, we will construct a cover of a $d$-ball of radius roughly $r/\wh{d}(a,\infty)^2$ by $d$-balls of radius roughly $r_i/\wh{d}(a,\infty)^2$.  Since the same factor appears in both the covered and covering balls, the resulting bounds on the Hausdorff measure of $\ovl{B}_{X,\wh{d}}(a,r)$ will be independent of $a$. 

An application of \eqref{inclusions hat middle} shows that
$$\ovl{B}_{d}\left( a, \frac{r}{\wh{d}(a,\infty)^2}\cdot \frac{C}{C+1} \right) \subeq \ovl{B}_{\wh{d}}(a,r).$$
We may also apply  \eqref{inclusions hat middle} to each $\wh{B}_i$, because
$$r_i < \frac{\wh{d}(a,\infty)-2r}{C} \leq \frac{\wh{d}(a,x_i)-2r + \wh{d}(x_i,\infty)}{C} \leq \frac{\wh{d}(x_i, \infty)}{C}.$$
Accordingly, 
$$\wh{B}_i \subeq \ovl{B}_{d}\left(x_i, \frac{r_i}{\wh{d}(x_i,\infty)^2}\cdot \frac{4C}{C-1}\right).$$
We also know that 
$$\wh{d}(x_i,\infty)^2 \geq (\wh{d}(a,\infty) - \wh{d}(a,x_i))^2 \geq (\wh{d}(a,\infty) - r - \epsilon)^2
\geq \wh{d}(a,\infty)^2\left(\frac{C-1}{C} - \frac{\epsilon}{\wh{d}(a,\infty)}\right)^2 >0,$$
and hence
$$\wh{B}_i \subeq B_i := \ovl{B}_{d}\left(x_i, \frac{r_i}{\wh{d}(a,\infty)^2}\cdot \frac{4C}{C-1}\left(\frac{C-1}{C} - \frac{\epsilon}{\wh{d}(a,\infty)}\right)^{-2}\right).$$
Note that the radius of each $B_i$ is bounded above by a constant which is independent of $i \in I$ and tends to zero as $\ep$ tends to zero.   Furthermore the collection $\{B_i\}_{i \in I}$ covers the ball
$$\ovl{B}_{d}\left(a,\frac{r}{\wh{d}(a,\infty)^2}\cdot \frac{C}{C+1}\right) .$$  
This implies that
$$\Hdim^Q_{d}\left(\ovl{B}_{d}\left(a, \frac{r}{\wh{d}(a, \infty)^2} \cdot \frac{C}{C+1} \right)\right)  \leq \frac{4C^3}{(C-1)^3\wh{d}(a, \infty)^2} \Hdim^Q_{\wh{d}}(\ovl{B}_{\wh{d}}(a, r)).$$
The Ahlfors $Q$-regularity of $(X,d)$ now yields 
\begin{equation}\label{Case 1 Lower}
\left(\frac{(C-1)^3}{4C^2(C+1)}\right)^Q\frac{r^Q}{K} \leq \mathcal{H}^{Q}_{\wh{d}}(\ovl{B}_{X,\wh{d}}(a,r)).
\end{equation}

To construct an upper bound for $\mathcal{H}^{Q}_{\wh{d}}(\ovl{B}_{\wh{d}}(a,r))$, we fix a new $\epsilon > 0$ such that 
$$ \epsilon < \min\left(\frac{1}{2\wh{d}(a,\infty)(\til{C}+1)}, \frac{1}{\til{C}+1}\right).$$ 
Consider a cover of $$\ovl{B}_d\left(a, \frac{r}{\wh{d}(a,\infty)^2} \cdot \frac{4\til{C}}{\til{C}-1}\right)$$ by a collection of balls $\{B_i\}_{i \in I}$ where for each $i \in I$, $B_i := \ovl{B}_d(x_i,r_i)$, $r_i < \epsilon$, and 
$$x_i \in B_d\left(a, \frac{r}{\wh{d}(a,\infty)^2} \cdot \frac{4\til{C}}{\til{C}-1} + \epsilon \right). $$
We will construct a cover of a $\wh{d}$-ball of radius $r$ by $\wh{d}$-balls of radius roughly $r_i\wh{d}(a,\infty)^2$.  This time, the cancellation will come from the fact that the $d$-ball we begin with has radius roughly $r/\wh{d}(a,\infty)^2$.

As $C > \til{C}$, we have $r \leq \wh{d}(a,\infty)/\til{C}$ as well, and so \eqref{inclusions hat middle} implies
$$\ovl{B}_{\wh{d}}(a,r) \subeq \ovl{B}_d \left( a, \frac{r}{\wh{d}(a,\infty)^2} \cdot \frac{4\til{C}}{\til{C}-1}\right).$$
We wish to apply \eqref{inclusions d middle} to each $B_i$, using the constant $\til{C}$.  Note that for every point $z \in X$, $\wh{d}(z,\infty)\leq 1$, and so the requirement that $r_i \leq (\wh{d}(x_i,\infty)(\til{C}+1))\inv$ is satisfied because $\ep < 1/(\til{C}+1)$.  Accordingly, for each $i \in I,$ 
\begin{equation} \label{not bi hat} B_i \subeq \ovl{B}_{\wh{d}}\left(x_i,r_i\wh{d}(x_i,\infty)^2 \cdot \frac{\til{C}+1}{\til{C}} \right).\end{equation}
We now wish to estimate $\wh{d}(x_i,\infty)^2$ independently of $i \in i$.  We have assumed that $r \leq \frac{\wh{d}(a,\infty)}{C}$, and so
\begin{equation}\label{r bound} r \leq \frac{\wh{d}(a,\infty)(\til{C}-1)}{8\til{C}(\til{C}+1)}.\end{equation}
The upper bound on $\ep$ and \eqref{r bound} are exactly what is needed to show that 
$$ \frac{r}{\wh{d}(a,\infty)^2} \cdot \frac{4\til{C}}{\til{C}-1} + \epsilon \leq \frac{1}{\wh{d}(a,\infty)(\til{C}+1)}.$$
Invoking \eqref{inclusions d middle}, we have
\begin{equation}\label{this thing}x_i \in B_d\left(a, \frac{r}{\wh{d}(a,\infty)^2} \cdot \frac{4\til{C}}{\til{C}-1} + \epsilon \right) \subeq B_{\wh{d}}\left(a, r \cdot \frac{4(\til{C} + 1)}{(\til{C}-1)} + \epsilon\wh{d}(a,\infty)^2\cdot \frac{(\til{C}+1)}{\til{C}}\right). \end{equation}
Now \eqref{r bound} and \eqref{this thing} provide the estimate
$$\wh{d}(x_i,\infty)^2 \leq (\wh{d}(a,x_i) + \wh{d}(a,\infty))^2  \leq \wh{d}(a,\infty)^2 \left(\frac{2\til{C}+1}{2\til{C}} + \frac{\epsilon \wh{d}(a,\infty)^2(\til{C} +1)}{\til{C}} \right)^2.  $$
Substituting this into \eqref{not bi hat}, we have
$$B_i \subeq \wh{B}_i:= \ovl{B}_{\wh{d}}\left(x_i, r_i \wh{d}(a,\infty)^2 \left(\frac{2\til{C}+1}{2\til{C}} + \frac{\epsilon \wh{d}(a,\infty)^2(\til{C} +1)}{\til{C}} \right)^2 \frac{(\til{C}+1)}{\til{C}}\right).$$
Note that the radius of each $\wh{B}_i$ is bounded above by a constant which is independent of $i \in I$ and tends to zero as $\ep$ tends to zero.  Moreover, the collection $\{\wh{B}_i\}_{i \in I}$ covers the ball $\ovl{B}_{\wh{d}}(a,r)$.  This combined with the Ahlfors $Q$-regularity of $(X,d)$ and \eqref{Case 1 Lower} shows that for $r \leq \wh{d}(a,\infty)/C$, 
\begin{equation}\label{Case 1}
\left(\frac{(C-1)^3}{4C^2(C+1)}\right)^Q\frac{r^Q}{K} \leq \mathcal{H}^{Q}_{\wh{d}}\left(\ovl{B}_{\wh{d}}(a,r)\right) \leq K \left(\frac{(2\til{C} +1)^2(\til{C} +1)}{\til{C}^2(\til{C}-1)}\right)^Q r^Q.
\end{equation}

In order to estimate $\Hdim^Q(\ovl{B}_{\wh{d}}(a,r))$ in the case that $r > \wh{d}(a,\infty)/C$, we first estimate the Hausdorff $Q$-measure of balls centered at $\infty$.   Consider $\ovl{B}_{\wh{X},\wh{d}}(\infty,r) \cap X$, where $r \leq \text{diam}(\wh{X},\wh{d})\leq 2.$   For $m \in \nats$, define the half-open annulus
$$A_m := B_{d}\left(p,\frac{2^{m+1}}{r}-1\right) - B_{d}\left(p,\frac{2^{m+1}}{r}-1\right).$$
By Remark \ref{infty balls}, we also have
$$A_m = \left(\ovl{B}_{\wh{X},\wh{d}}\left(\infty, \frac{r}{2^m}\right) - \ovl{B}_{\wh{X},\wh{d}}\left(\infty, \frac{r}{2^{m+1}}\right)\right) \cap X,$$
showing that $\ovl{B}_{\wh{X},\wh{d}}(\infty,r) \cap X = \bigcup_{m=0}^\infty A_m.$  As the collection $\{A_m\}_{m \in \nats}$ is disjointed, it is sufficient to find suitable bounds on the $Q$-Hausdorff measure of each $A_m$.  Lemma \ref{inclusions} and a covering argument similar to those above yield that 
$$ \mathcal{H}_{\wh{d}}^Q(A_m) \leq K\left(\frac{\til{C}+1}{\til{C}}\right)^Q \left(\frac{r}{2^{m-1}}\right)^Q.$$
Summing over $m \in \nats$, we may conclude that 
\begin{equation}\label{Case 2 Upper}
\mathcal{H}^Q_{\wh{X},\wh{d}}(\ovl{B}_{\wh{X}, \wh{d}}(\infty,r)\cap X) \leq K\left(\frac{4(\til{C}+1)}{\til{C}}\right)^Q \frac{r^Q}{2^Q-1}.
\end{equation}
Note that \eqref{Case 2 Upper} also holds for $r \geq \diam{(X,\wh{d})}.$

To establish a lower bound, we first note that $X$ is connected, $\wh{d}(p,\infty)=1$, and $r \leq 2$, and therefore there exists a point $a \in X$ such that $\wh{d}(a,p)=1-r/2.$  Lemma \ref{p symmetry 1} shows that $\wh{d}(a,\infty) = r/2$.  Noting that 
$$\ovl{B}_{\wh{d}}\left(a,r/(2C)\right) \subeq (\ovl{B}_{\wh{X},\wh{d}}(\infty,r) \cap X),$$
applying \eqref{Case 1} to $\ovl{B}_{\wh{d}}(a, \frac{r}{2C})$ shows that 
\begin{equation}\label{Case 2}
\left(\frac{(C-1)^3}{8C^3(C+1)}\right)^Q\frac{r^Q}{K} \leq \mathcal{H}^Q_{\wh{X}, \wh{d}}\left(\ovl{B}_{\wh{X}, \wh{d}}(\infty,r)\right) \cap X.
\end{equation}

Finally, we consider $\ovl{B}_{\wh{d}}(a,r)$ with $a \in X$ and $r>\wh{d}(a,\infty)/C.$  In this case, 
$$\ovl{B}_{\wh{d}}(a,r)  \subeq \ovl{B}_{\wh{X},\wh{d}}\left(\infty,(C+1)r\right) \cap X,$$
and so \eqref{Case 2 Upper} provides the upper bound
\begin{equation}\label{Case3Upper}
\mathcal{H}^Q_{\wh{d}}\left(\ovl{B}_{\wh{d}}(a,r)\right) \leq  K\left(\frac{4(\til{C}+1)^2}{\til{C}}\right)^Q \frac{r^Q}{2^Q-1}.
\end{equation}
For the lower bound, we first suppose that $r \leq C\wh{d}(a,\infty).$  Applying \eqref{Case 1 Lower} to the $\ovl{B}_{\wh{d}}(a, r/C^2)$ shows that 
$$\left(\frac{(C-1)^3}{4C^4(C+1)}\right)^Q\frac{r^Q}{K} \leq \mathcal{H}^Q_{\wh{d}}\left(\ovl{B}_{\wh{d}}(a,r)\right).$$
If instead $r > C\wh{d}(a,\infty),$ then 
$$\ovl{B}_{\wh{d}}\left(\infty,\frac{C-1}{C}r \right) \cap X \subeq \ovl{B}_{\wh{d}}(a,r).$$
Thus we may apply \eqref{Case 2} above to show that 
$$\left(\frac{(C-1)^4}{8C^4(C+1)}\right)^Q\frac{r^Q}{K} \leq \mathcal{H}^Q_{\wh{d}}\left(\ovl{B}_{\wh{d}}(a,r)\right).$$
In either case,
\begin{equation}\label{Case 3}
\left(\frac{(C-1)^3}{8C^4(C+1)}\right)^Q\frac{r^Q}{K} \leq \mathcal{H}^Q_{\wh{d}}\left(\ovl{B}_{\wh{d}}(a,r)\right) \leq K\left(\frac{4(\til{C}+1)^2}{\til{C}}\right)^Q \frac{1}{2^Q-1}r^Q.
\end{equation}

Estimates \eqref{Case 1} and \eqref{Case 3} show that $(X, \wh{d})$ is Ahlfors $Q$-regular; different choices of $\til{C}$ will yield different constants.  However, any choice of $\til{C}$ yields a constant which is $O(K)$ for a fixed dimension $Q$.   
\end{proof}

\begin{prop}\label{hat LLC}  Let $(X,d)$ be an unbounded and $\lambda$-$LLC$ metric space.  Then $(X,\wh{d})$ is $\lambda'$-$LLC$, where $\lambda'$ depends only on $\lambda$.  
\end{prop}

\begin{proof}  This follows immediately from Lemma \ref{QM} and Theorem \ref{LLC is qs invariant}.  It can also be shown directly using Lemma \ref{inclusions}.  
\end{proof}

\begin{proof}[Proof of Theorem \ref{main} (iv) and (v)]
Let $(X,d)$ be an Ahlfors $2$-regular, $LLC$ metric space which is homeomorphic to $\reals^2$.   Then by Remark \ref{infty balls} and Propositions \ref{hat regular} and \ref{hat LLC},  $(X,\wh{d})$ is a bounded, Ahlfors $2$-regular and $LLC$ metric space with data depending only on the data of $(X,d)$.  

If $(X,d)$ is complete, Lemma \ref{hat boundary} shows that $\card{\wh{\partial}{X}} = 1$.   Thus by Theorem \ref{main} (ii), there is a $\eta$-quasisymmetric homeomorphism $f\colon (X,\wh{d})\to \sphere^*$ where $\eta$ depends only on the data of $X$.  By Lemma \ref{QM} there is a quasi-M\"obius homeomorphism $\iota\colon (X,d) \to (X,\wh{d})$.  As the quasi-M\"obius condition is invariant under composition with bi-Lipschitz maps, Lemmas \ref{QM} and \ref{plane to sphere} show that there is a quasi-M\"obius homeomorphism  $g\colon \sphere^* \to \reals^2$.   The composition $g\circ f \circ \iota$ sends unbounded sequences to unbounded sequences, and so Theorem \ref{QM-QS} shows that $g \circ f \circ \iota$ is quasisymmetric.  As $\iota$ and $g$ have fixed distortion functions, the distortion function of $g \circ f \circ \iota$ depends only on $\eta$ and hence only on the data of $X$.  

Now suppose that $(X,d)$ is not complete.  Lemma \ref{hat boundary} shows that $\card{\wh{\partial}X} \geq 2$.  In the construction of the warped space $(X,\wh{d})$, the basepoint $p$ can be chosen arbitrarily; as $(\ovl{X},d)$ is connected and unbounded, we may choose $p$ such that there is a point $z \in \partial{X}$ such that $d(z,p) \leq 1$.  Then $\wh{d}(\infty,z)\geq 1/2$, and so by Remark \ref{infty balls} and Lemma \ref{hat boundary} we have
$$ 1 \leq \frac{\diam{\wh{X}}}{\diam{\wh{\partial}{X}}} \leq 4.$$  
By Theorem \ref{main} (iii), there is an $\eta$-quasisymmetric homeomorphism $f\colon (X,\wh{d}) \to \disk$ where $\eta$ depends only the data of $X$ and the ratio $\diam{\wh{X}}/\diam{\wh{\partial}{X}}$.  Since this ratio is bounded above and below, $\eta$ depends only on the data of $X$.   We may choose a M\"obius homeomorphism $g\colon \disk \to \half$ so that the composition $g \circ f \circ \iota$ maps unbounded sequences to unbounded sequences.  As in the complete case, it follows that $g \circ f \circ \iota$ is $\eta'$-quasisymmetric where $\eta'$ depends only on $\eta$.  
\end{proof}

\bibliographystyle{amsplain}
\bibliography{QSbibli}

\end{document}